\documentclass[11pt]{amsart}

\usepackage{amsthm}
\usepackage{amsmath}
\usepackage{amsfonts}
\usepackage{amssymb}

\usepackage{fullpage}

\newcommand{\nc}[2]{\newcommand{#1}{#2}}

\newcommand{\ie}{i.e.\ }

\newcommand{\cf}{cf.\ }

\nc{\biggiven}{\; \bigg|\,}
\nc{\biggergiven}{\; \Bigg|\;}
\nc{\given}{\,|\,}
\nc{\bigggiven}{\; \right|\left.\,}
\nc{\bgiven}{\; \big| \,}

\newcommand{\ba}{\begin{align*}}
\newcommand{\ea}{\end{align*}}
\def\ba #1\ea{\begin{align*}#1\end{align*}}
\newcommand{\inv}[1]{#1^{-1}}
\nc{\indfn}{\mathbf{1}}
\newcommand{\ind}[1]{\indfn_{\left\{#1\right\}}}
\newcommand{\cl}[1]{\overline{#1}} % closure operator
\newcommand{\rvect}[1]{\boldsymbol #1}
\newcommand{\vect}[1]{\rvect{#1}}
\newcommand{\abs}[1]{\lvert #1 \rvert}

\newcommand{\comp}[1]{#1^{\mathrm{c}}}

\renewcommand{\P}{\mathsf{P}} % underlying prob measure
\newcommand{\EP}{\mathsf{E}} % expectation operator wrt \P

\nc{\cont}{\subset}
\nc{\goesto}{\rightarrow}
\nc{\conv}{\longrightarrow}
\nc{\wc}{\Rightarrow}
\nc{\vconv}{\stackrel{v}{\longrightarrow}}
\nc{\vgoesto}{\stackrel{v}{\goesto}}
\newcommand{\mt}[1]{\ \ \quad\text{#1}\ \,}
\newcommand{\mtt}[1]{\ \ \mathit{#1}\ }
\nc{\RV}{\text{RV}}
\nc{\eqdist}{\stackrel{d}{=}}
\nc{\iid}{\stackrel{\text{iid}}{\sim}}
\nc{\bdry}{\partial}
\nc{\mand}{\qquad\text{and}\qquad}
\nc{\clos}{\text{cl}}
\nc{\smax}{\vee}
\nc{\smin}{\wedge}

\newcommand{\pp}[1]{\epsilon_{#1}}

\nc{\R}{\mathbb{R}}
\newcommand{\E}{\mathbb{E}}
\newcommand{\mplus}[1]{\mathbb{M}_{+}#1}
\nc{\Bor}{\mathcal{B}}
\nc{\contfn}{\mathcal{C}^+_K}
\nc{\cpt}{\mathcal{K}}
\nc{\cbfn}{\mathcal{C}}

\theoremstyle{plain}
\newtheorem{thm}{Theorem}[section]
\newtheorem{lem}{Lemma}[section]
\newtheorem{cor}{Corollary}[section]
\newtheorem{prop}{Proposition}[section]

\theoremstyle{definition}
\newtheorem*{defn}{Definition}
\newtheorem*{rmk}{Remark}

\newenvironment{pf}[1][Proof]{\begin{proof}[\textit{\textbf{#1}}]} {\end{proof}}

\newtheoremstyle{example}{10pt}{3pt}% Spacing
	{}% Body font
	{}% Indent
	{\itshape\bfseries}% Thm head font
	{.}%        Punctuation after thm head
	{0.5em}%     Space after thm head
	{}%         Thm head spec (can be left empty, meaning `normal')

\theoremstyle{example}
\newtheorem*{exnotnumbered}{Example}
\newtheorem{exnumbered}{Example}[section]
\newenvironment{exnum}{\pushQED{\qed}\begin{exnumbered}} {\qedhere\end{exnumbered}}

\usepackage{ifpdf}
\ifpdf
\usepackage{hyperref}
\hypersetup{%pdfpagemode=FullScreen,%
	pdftitle={Asymptotics of Markov Kernels and the Tail Chain},%
	pdfauthor={S. I. Resnick and D. Zeber},%
	pdfcreator={},%
	pdfproducer={},}
\fi

\newcommand{\sid}[1]{#1}

\numberwithin{equation}{section}
\numberwithin{figure}{section}

\nc{\bdist}{H}
\nc{\Gz}{G(\{0\})}
\nc{\eb}{y}
\nc{\ebt}{\eb(t)}
\nc{\Xv}{\rvect{X}}
\nc{\xv}{\rvect{x}}
\nc{\bx}{\rvect{x}}
\nc{\bX}{\rvect{X}}
\nc{\bz}{\rvect{0}}
\nc{\binfty}{\rvect{\infty}}
\nc{\bT}{\rvect{T}}
\newcommand{\mckern}[3][K]{#1\big(#2\,,\, #3\big)} 
\newcommand{\mckernl}[3][K]{#1(#2\,,\, #3)}

\title{Asymptotics of Markov Kernels and the Tail Chain}
\author[S.~I.~Resnick]{Sidney I.~Resnick}
\author[D.~Zeber]{David Zeber}

\address{Sidney I.~Resnick, School of ORIE, Rhodes Hall 284, Cornell University,  Ithaca, NY 14853}
\email{sir1@cornell.edu}%
\address{David Zeber, Department of Statistical Science, 301 Malott Hall, Cornell University, Ithaca NY 14853}
\email{dsz5@cornell.edu}%

\thanks{S.~I.~Resnick and D.~Zeber were partially supported by ARO   Contract W911NF-10-1-0289 and NSA Grant H98230-11-1-0193 at Cornell University.}
 
\date{\today}

\begin{document}

\begin{abstract}
An asymptotic model for extreme behavior of certain Markov chains is the ``tail chain''.
Generally taking the form of a multiplicative random walk, it is useful in deriving extremal characteristics such as point process limits. 
We place this model in a more general context, formulated in terms of extreme value theory for transition kernels, and 
extend it 
by formalizing the distinction between extreme and non-extreme states. 
We make the link between the update function and transition kernel forms considered in previous work, and we show that the tail chain model 
leads to a multivariate regular variation property of the finite-dimensional distributions under assumptions on the marginal tails alone. 
\end{abstract}

%%%%%%%%%%%%%%%%%%%%%%%%%%%%%%%%%%%%%%

\maketitle

%%%%%%%%%%%%%%%%%%%%%%%%%%%%%%%%%%%%%%

\section{Introduction}

A method of approximating the extremal behavior of discrete-time Markov chains is to use an asymptotic process called the \emph{tail chain} under an asymptotic assumption on the transition kernel of the chain.
Loosely speaking, if the distribution of the next state converges
under some normalization as the current state becomes extreme, then
the Markov chain behaves approximately as a multiplicative random walk
upon leaving a large initial state. This approach leads to intuitive
extremal models in such cases as autoregressive processes with random
coefficients, which include a class of ARCH models.  
The focus on Markov kernels was introduced by Smith \cite{smith1992extremal}. 
Perfekt \cite{perfekt1994extremal, perfekt1997extreme} extended the approach to higher dimensions, and Segers \cite{segers2007multivariate} rephrased the conditions in terms of update functions. 

Though  not  restrictive in practice, the previous approach
 tends to mask  aspects of the processes'
extremal behaviour.  
Markov chains which admit the tail chain approximation fall into one
of two  categories.  
Starting from an extreme state, the chain either remains extreme over any finite time horizon, or will drop to a ``non-extreme'' state of lower order after a finite amount of time. 
The latter case is problematic in that the tail chain model is not sensitive to possible subsequent jumps from a non-extreme state to an extreme one. 
% which are entirely possible. 
Previous developments handle this by ruling out the class of processes exhibiting this behaviour via a technical condition, which we refer to as the \emph{regularity condition}. 
Also, most previous work has assumed stationarity, since interest focused on computing the extremal index or deriving limits for the exceedance point processes,  drawing on the theory established for stationary processes with mixing by Leadbetter et al.\ \cite{leadbetter1983extremes}. 
However, stationarity is not fundamental in determining the extremal behaviour of the finite-dimensional distributions.

We place the tail chain approximation in the context of an extreme value theory for Markovian transition kernels, which \emph{a priori} does not necessitate any such restrictions on the class of processes to which it may be applied. 
In particular, we introduce the concept of boundary distribution, which controls tail chain transitions from non-extreme to extreme. 
Although distributional convergence results are more naturally phrased in terms of transition kernels, we treat the equivalent update function forms as an integral component to interfacing with applications, and we phrase relevant assumptions in terms of both. 
While not making explicit a complete tail chain model for the class of chains excluded previously, we demonstrate the extent to which previous models may be viewed as a partial approximation within our framework. 
This is accomplished by formalizing the division between extreme and non-extreme states as a level we term the  \emph{extremal boundary}. 
We show that, in general, the tail chain approximates the \emph{extremal component}, the portion of the original chain having yet to cross below this boundary. 
Phrased in these terms, the regularity condition requires that the distinction between the original chain and its extremal component disappears asymptotically.

After introducing our extreme value theory for transition kernels, along with a representation in terms of update functions, we derive limits of finite-dimensional distributions conditional on the initial state, as it becomes extreme. 
We then examine the effect of the regularity condition on these results. Finally, adding the assumption of marginal regularly varying tails leads to convergence results for the unconditional distributions akin to regular variation.

%%%%%%%%%%%%%%%%%%%%%%%%%%%%%%%%%%%%%%%%%%%%%%%%

\subsection{Notation and Conventions}
We review notation and  relevant concepts. If not explicitly
specified, assume that any space $\mathbb{S}$ under discussion is a
topological space paired with its Borel $\sigma$-field of open sets
$\Bor(\mathbb{S})$ to form a measurable space.   
Denote by $\cpt(\mathbb{S})$ the collection of its compact sets; by $\cbfn(\mathbb{S})$ the space of real-valued continuous, bounded functions on $\mathbb{S}$; and by $\contfn(\mathbb{S})$ the space of non-negative continuous functions with compact support. Weak convergence of probability measures is represented by $\wc$.

For a space $\E$ which is locally compact with countable base (for
example, a subset of $[-\infty,\infty]^d$), $\mplus(\E)$  is  the
space of non-negative Radon measures on $\Bor(\E)$; point \sid{measures consisting of
single point masses at $x$}  will be written as $\pp{x}(\cdot)$. 
A sequence of measures $\{ \mu_n\} \cont \mplus(\E)$  converges
vaguely to $\mu\in \mplus(\E)$ (written $\mu_n\vgoesto \mu$) if
$\int_\E f\, d\mu_n \goesto \int_\E f\, d\mu$ as $n\goesto\infty$ for
any $f\in\contfn(\E)$.  
The shorthand $\mu(f) = \int f\; d\mu$ \sid{is} handy. 
That the  distribution of a random vector $\bX$ is regularly varying on a
cone $\E \cont [-\infty,\infty]^d \backslash \{\bz\}$
means that $t\,\P[\bX/b(t) \in \cdot\, ] \vgoesto
\mu^*(\cdot)$ in $\mplus(\E)$ as $t\goesto\infty$ for some
non-degenerate limit measure $\mu^*\in \mplus(\E)$ and scaling
function $b(t) \goesto \infty$. The limit $\mu^*$ is necessarily homogeneous in the sense that
$\mu^*(c\,\cdot) = c^{-\alpha}\mu^*(\cdot)$ for some $\alpha>0$.
The regular variation is \emph{standard} if $b(t)=t$.

If $\bX = (X_0,X_1,X_2,\dots)$ is a (homogeneous) Markov chain and $K$
is a Markov transition kernel, we write $\bX \sim K$ to mean that the
dependence structure of $\bX$ is specified by $K$, \ie 
\[ \P[X_{n+1} \in \cdot \given X_{n} = x ] = \mckern{x}{\cdot},\qquad
n=0,1,\dotsc .
\] 
We adopt the standard shorthand $\P_x [(X_1,\dots,X_m) \in \cdot\,] =
\P[(X_1,\dots,X_m) \in \cdot \given X_0 = x]$.
Some useful technical results are assembled in Section \ref{seclem} (p.~\pageref{seclem}).

%%%%%%%%%%%%%%%%%%%%%%%%%%%%%%%%%%%%%%

\section{Extremal Theory for Markov Kernels} 
We begin by focusing on the Markov transition kernels
 rather than the stochastic processes they determine,  
and introduce a class of kernels we
term ``tail kernels,''  
which we will view as \sid{scaling} limits of certain kernels. 
Antecedents include Segers'  \cite{segers2007multivariate} 
definition of ``back-and-forth tail chains'' that 
approximate certain Markov chains started from an extreme
value. 

\sid{For} a Markov chain $\bX\sim K$ on $[0,\infty)$, it is reasonable
to expect that extremal behaviour of $\bX$ is determined by pairs
 $(X_n,X_{n+1})$, and one way to control such pairs is to assume that
$(X_n,X_{n+1})$ belongs to a bivariate domain of attraction (\cf \cite{bortot2000sufficiency, smith1992extremal}). 
In the context of regular variation, 
writing
\begin{equation} \label{eqmchextr}
t\,\P\left[\frac{X_n}{b(t)}\in A_0\,,\, \frac{X_{n+1}}{b(t)} \in A_1 \right] = \int_{A_0} \mckern{b(t)u}{b(t)A_1} \; t\,\P\left[\frac{X_n}{b(t)} \in du \right]
\end{equation}
%which 
suggests combining marginal regular variation of $X_n$ with a scaling
kernel limit to derive extremal properties of the finite-dimensional
distributions 
 (fdds) 
\cite{perfekt1994extremal,
  perfekt1997extreme, segers2007multivariate}, and this is the direction
we take. We first discuss the kernel scaling  operation.

For simplicity, we  assume the state space of the Markov chain is
$[0,\infty)$, although with suitable modifications, it is relatively
straightforward to extend the results to $\R^d$.  
Henceforth $G$ and $H$ will denote probability distributions on $[0,\infty)$.

%%%%%%%%%%%%%%%%%%%%%%%%%%%%%%%%%%%%%%%%%%%%%%%%%%

\subsection{Tail Kernels}\label{subsec:tailkernels}
 
The \emph{tail kernel associated with} $G$, with \emph{boundary distribution} $\bdist$, is
\begin{equation} \label{eqtailkerndef}
 \mckern[K^*]{y}{A} = \begin{cases} G(\inv{y} A) & \ 
y>0 \\ \bdist(A) & \ 
y=0 \end{cases} 
\end{equation}
for any measurable set $A$. 
Thus, the class of tail kernels on $[0,\infty)$ is parameterized by
the pair of probability distributions $(G,\bdist)$. Such kernels are
characterized by a scaling property:

\begin{prop} \label{proptailkernhom}
A Markov transition kernel $K$ is a tail kernel associated with some $(G,\bdist)$ if and only if it satisfies the relation 
\begin{equation} \label{eqtailkernhom}
 \mckern{uy}{A} = \mckern{y}{\inv{u} A}
\end{equation}
when $y>0$ for any $u>0$, in which case $G(\cdot) = \mckernl{1}{\cdot}$. The property \eqref{eqtailkernhom} extends to $y=0$ iff $\bdist=\pp{0}$. 
\end{prop} 

\begin{pf}
If $K$ is a tail kernel, \eqref{eqtailkernhom} follows directly from the definition. Conversely, assuming \eqref{eqtailkernhom}, for $y>0$ we can write \[ \mckern{y}{A} = \mckern{1}{\inv{y} A}, \] demonstrating that $K$ is a tail kernel associated with $\mckernl{1}{\cdot}$ (with boundary distribution $\bdist = \mckernl{0}{\cdot}$). 
To verify the second assertion, fixing $u>0$, we must show that $\bdist(\inv{u}\cdot) = \bdist(\cdot)$ iff $\bdist=\pp{0}$. On the one hand, 
we have $\pp{0}(\inv{u}A) = \pp{0}(A)$. On the other, $\bdist(0,\infty) = \lim_{n\goesto\infty} \bdist(\inv{n},\infty) = \bdist(1,\infty)$, so $\bdist(0,1] = 0$. A similar argument shows that $\bdist(1,\infty)=0$ as well. 
\end{pf}

We call the Markov chain $\bT \sim K^*$  a \emph{tail chain associated
  with} $(G,H)$. Such a chain can be represented as
\begin{equation} \label{eqtailkernform}
T_n = \xi_n\, T_{n-1} + \xi'_n\, \ind{T_{n-1}=0} \quad\mt{for} n=1,2,\dotsc, 
\end{equation}
where $\xi_n \iid G$ and $\xi'_n \iid \bdist$ are independent of each other and of $T_0$. If $H=\pp{0}$, then $\bT$ becomes a multiplicative random walk with step distribution $G$ and absorbing barrier at $\{0\}$:  $T_n = T_0 \, \xi_1 \dotsm \xi_n$.

%%%%%%%%%%%%%%%%%%%%%%%%%%%%%%%%%%%%%%%%%%%%%%%%

\subsection{Convergence to Tail Kernels}

The tail chain approximates the behaviour of a Markov chain
$\bX\sim K$ in extreme states. Asymptotic results require 
that the normalized distribution of $X_1$ be well-approximated by some
distribution $G$ when $X_0$ is large, and we
interpret this requirement as a domain of attraction condition for kernels. 

\begin{defn}
A Markov transition kernel $K : [0,\infty) \times \Bor[0,\infty)\goesto [0,1]$ is in the \emph{domain of attraction of} $G$, written $K\in D(G)$, if as $t\goesto\infty$,
\begin{equation} \label{eqkerndoa}
 \mckern{t}{t\cdot} \wc G(\cdot) \mt{on}  [0,\infty] .
\end{equation} 
\end{defn}
\noindent Note that $D(G)$ contains at least the class of tail kernels
associated with $G$ (\ie with any boundary distribution $\bdist$).
A simple scaling argument extends \eqref{eqkerndoa} to 
\begin{equation} \label{eqkerndoascaling}
\mckern{tu}{t\cdot} \wc G(\inv{u} \cdot) =:\mckern[K^*]{u}{\cdot\,},\qquad  u>0,
\end{equation}
 where $K^*$ is any tail kernel associated with $G$; 
this is the form appearing in
 \eqref{eqmchextr}.  
 \sid{Thus tail} kernels
\sid{are scaling}
limits for kernels in a domain of attraction. 
In fact, tail kernels are the only possible limits: 

\begin{prop}
Let $K$ be a transition kernel and $\bdist$ be an arbitrary
distribution on $[0,\infty)$. If for each $u>0$ there exists a
distribution $G_u$ such that 
$ \mckern{tu}{t\cdot} \wc G_u(\cdot) $
as $t\goesto\infty$, then the function $\widehat{K}$ defined on
$[0,\infty)\times \Bor[0,\infty)$ as  
\[
 \mckern[\widehat{K}]{u}{A}
:= \begin{cases} G_u(A) & 
u>0 \\ \bdist(A) & 
u=0 \end{cases} \]  is a
tail kernel associated with $G_1$.  
\end{prop}

\begin{pf}
It suffices to show that $G_u(\cdot) = G_1(\inv{u}\cdot)$ for any $u>0$. 
But this follows directly from the uniqueness of weak limits, since \eqref{eqkerndoascaling} shows that $\mckernl{tu}{t\,\cdot} \wc G_1(\inv{u}\cdot)$. 
\end{pf}

\noindent A version of \eqref{eqkerndoascaling} 
uniform in $u$ is needed for 
 fdd convergence results. 
\begin{prop} \label{propkernconvunif}
Suppose $K \in D(G)$, and $K^*$ is a tail kernel associated with $G$. Then, for any $u>0$ and any non-negative function $u_t = u(t)$ such that $u_t \goesto u$ as $t\goesto\infty$, we have 
\begin{equation} \label{eqkernunifconv}
\mckern{tu_t}{t\cdot} \wc \mckern[K^*]{u}{\cdot\,} ,\qquad
(t\goesto\infty). 
\end{equation}
\end{prop}

\begin{pf}
Suppose $u_t\goesto u>0$. Observe that $\mckernl{tu_t}{t\,\cdot} =
\mckernl{tu_t}{(tu_t)\, \inv{u_t} \cdot}$, and put $h_t(x) = u_tx$,
$h(x) = ux$. Writing $P_t(\cdot) = \mckernl{tu_t}{tu_t\,\cdot}$, we
have
\[ \mckern{tu_t}{t\cdot} = P_t \circ \inv{h_t} \wc G \circ
\inv{h} = G(\inv{u}\cdot) = \mckern[K^*]{u}{\cdot\,} 
\hfil
\]   by \cite[Theorem 5.5, p.~34]{billingsley1968convergence}. 
\end{pf}

\sid{The measure} $G$ controls $\bX$ upon leaving an extreme state,
and $H$  describes the possibility of jumping from a non-extreme
state to an extreme one.  
The traditional assumption \eqref{eqkerndoa}  provides no
information 
\sid{about}
$\bdist$, and in fact \eqref{eqkernunifconv} may fail if $u=0$---see
Example \ref{exnonunif}. However, the choice of $H$ \sid{cannot  be
ignored if $0$ is an accessible point of the state space,} especially
for cases where $\Gz =
\mckernl[K^*]{y}{\{0\}} > 0$. 
We propose pursuing implications of the traditional
assumption \eqref{eqkerndoa} alone, and will add conditions as needed
to understand boundary behaviour of $\bX$.

Alternative, more general formulations of \eqref{eqkerndoa} include
replacing $\mckernl{t}{t\,\cdot}$ with  $\mckernl{t}{a(t)\, \cdot}$ or
$\mckernl{t}{a(t) \,\cdot + \, b(t)}$ with appropriate  functions
$a(t)>0$ and $b(t)$, in analogy with the \sid{usual} domains of attraction
conditions in extreme value theory. 
Indeed, the second choice coincides with the original presentation by
Perfekt \cite{perfekt1994extremal}, and relates to the conditional
extreme value model  \cite{das2011conditioning,
  heffernan2007limit,
heffernan2004conditional}. For clarity, and to maintain ties with
regular variation, we retain
the standard choice $a(t)= t$, $b(t)=0$.

%%%%%%%%%%%%%%%%%%%%%%%%%%%%%%%%%%%%%%%%%%%%%%%%%%%%%%%%%%

\subsection{Representation}
\sid{How do we} characterize  kernels belonging to
$D(G)$?
{}From \eqref{eqtailkernform}, for chains transitioning according to a
tail kernel, the next state is a random multiple of
the previous one, \sid{provided} the prior state is non-zero.  
We expect that chains transitioning according to 
$K\in D(G)$ behave approximately like this
upon leaving a large state, and this is best expressed 
in terms of a function describing how a new state
depends on the prior one.

Given a kernel $K$, we can always find a sample space $\E$, a
measurable  function $\psi : [0,\infty) \times \E \goesto
[0,\infty)$ and an $\E$-valued random element $V$ such that 
$ \psi(y, V) \sim \mckernl{y}{\cdot\,}$ for all $y$. 
Given a random variable $X_0$, 
if we define the process $\bX = (X_0,X_1,X_2,\dots)$ recursively as 
\[ X_{n+1} = \psi(X_n, V_{n+1}), \qquad n\geq 0,
\] 
where $\{V_n\}$ is an iid sequence equal in distribution to $V$ and independent of $X_0$, then $\bX$ is a Markov chain with transition kernel $K$. 
Call the function $\psi$ an \emph{update function corresponding to} $K$.  If in addition $K\in D(G)$, the domain of attraction condition
\eqref{eqkerndoa} 
becomes 
\[ \inv{t} \psi(t, V) \wc \xi, \] where $\xi\sim G$. 
Applying 
the probability integral transform or the Skorohod representation theorems \cite[Theorem 3.2, p.~6]{billingsley:1971}, \cite[Theorem 6.7, p.~70]{billingsley1999convergence}, we get the following result.

\nc{\newpm}{\lambda}

\begin{prop} \label{propupdfnas}
If $K$ is a transition kernel,  $K\in D(G)$ if and only if
there exists a measurable function $\psi^*:[0,\infty)\times [0,1] \goesto [0,\infty)$ and a random variable $\xi^*\sim G$ on the uniform
probability space $([0,1],\Bor,\newpm)$ such that   
\begin{equation} \label{equpdfnas}
t^{-1}{\psi^*(t,u)} 
\conv
 \xi^*(u) 
\qquad \forall \,u \in [0,1]
\end{equation}
as $t\goesto\infty$, and $\psi^*$
is an update function corresponding to $K$ in the sense that \[ \newpm\big[ \psi^*(y,\cdot)\in A\, \big] =
\mckern{y}{A}
\] for measurable sets $A$. 
\end{prop}

\noindent Think of the update function as $\psi^*(y,U)$ where $U(u)= u$ is a
uniform random variable on $[0,1]$.

\begin{pf}
If there exist such $\psi^*$ and $\xi^*$ satisfying \eqref{equpdfnas} then clearly $K\in D(G)$. 
Conversely, suppose $\psi(\cdot,V)$ is an update function corresponding to $K$. 
According to Skorohod's representation theorem (\cf Billingsley  \cite{billingsley1999convergence} p.\ 70, with the necessary modifications to allow for an uncountable index set), there exists a random variable $\xi^*$ and a stochastic process $\{ Y^*_t \,;\, t\geq 0\}$ defined on  
the uniform probability space 
$([0,1],\Bor,\newpm)$, 
taking values in $[0,\infty)$, 
such that 
\[ \xi^*\sim G\,, \qquad  Y^*_0 \eqdist \psi(0,V)\,,\qquad   Y^*_t \eqdist \inv{t}\psi(t,V)\ \ \text{for } t>0, \]
and $Y^*_t(u) \goesto \xi^*(u)$ as $t\goesto\infty$ for every $u\in[0,1]$. 
Now, define $\psi^* : [0,\infty)\times [0,1] \goesto [0,\infty)$ as 
\[ \psi^*(0,u) = Y^*_0(u) \qquad\text{and}\qquad \psi^*(t,u) = tY^*_t(u)\,, \ \  t> 0, \qquad \forall\, u\in [0,1]. \] 
It is evident that $\newpm[\psi^*(y,\cdot) \in A] = \P[\psi(y,V)\in A]$ for $y\in [0,\infty)$, so $\psi^*$ is indeed an update function corresponding to $K$, and $\psi^*$ satisfies \eqref{equpdfnas} by construction.  
\end{pf}

Update functions corresponding to $K$ are
not unique, and  some of them may fail
to converge pointwise as in \eqref{equpdfnas}. However
\eqref{equpdfnas} is convenient, and Proposition \ref{propupdfnas} shows that Segers' \cite{segers2007multivariate}
Condition 2.2  in terms of update
functions is equivalent to our weak convergence formulation  $K\in
D(G)$.

Pointwise convergence in \eqref{equpdfnas} gives an intuitive representation of kernels in a domain of attraction.

\begin{cor} \label{corupdfnrep}
$K\in D(G)$ iff
there exists a random variable $\xi\sim G$ defined on the uniform probability space, and a measurable function $\phi : [0,\infty)\times [0,1] \goesto (-\infty,\infty)$ satisfying $\inv{t}\phi(t,u) \goesto 0$ for all $u\in[0,1]$
such that   
\begin{equation} \label{equpdfnrep}
 \psi(y,u) 
 := \xi(u)\, y + \phi(y,u) 
\end{equation}
is an update function corresponding to $K$. 
\end{cor}

\begin{pf}
If such $\xi$ and $\phi$ exist, then $\inv{t} \psi(t,u) = \xi(u) + \inv{t}\phi(t,u) \goesto \xi(u)$ for all $u$, so $\psi$ satisfies \eqref{equpdfnas}. 
The converse follows from \eqref{equpdfnas}.
\end{pf}

Many Markov chains such as ARCH, GARCH and autoregressive processes
are specified by structured recursions that allow quick recognition of
update functions  corresponding to kernels in a domain of attraction.
A common example is the update function 
$ \psi(y, (Z,W)) = Zy + W$, which behaves like $\psi'(y,Z) = Zy$ when $y$ is large---compare $\psi'$ to the form \eqref{eqtailkernform} discussed for tail kernels. 
In general, if $K$ has an update function $\psi$ of the form 
\begin{equation} \label{equpdstdform}
 \psi(y, (Z,W)) = Zy + \phi(y,W) 
\end{equation}
for a random variable $Z\geq 0$ and a random element $W$, where
$\inv{t}\phi(t,w)\goesto 0$ whenever $w\in C$ for which $\P[W\in C]=1$, then $K\in D(G)$ with $G = \P[Z \in \cdot\,]$. We will refer to update functions 
satisfying \eqref{equpdstdform} as being in \emph{canonical form}.

%%%%%%%%%%%%%%%%%%%%%%%%%%%%%%%%%%%%%%%%%%%%%%%%%%%%%%%

\section{Finite-Dimensional Convergence and the Extremal Component}
\label{secfdd}

Given a Markov chain $\bX \sim K \in D(G)$,
we show that the finite-dimensional distributions 
(fdds) of $\bX$, started
from an extreme state,  converge to those of the tail chain $\bT$
defined in \eqref{eqtailkernform}. We initially develop results that
depend only
on $G$  (but not $H$), and then clarify what behaviour of $\bX$ is
controlled by  $G$ and $H$ respectively. 
We  make explicit links with prior work that did not consider the notion of 
boundary distribution.

If $\Gz=0$, the choice of $\bdist$ is inconsequential, since
$\P\big[\bT \text{ eventually hits } \{0\} \big] = 0$ and 
 $\bT$ is indistinguishable from the multiplicative random walk
 $\{T^*_n = T_0 \xi_1\cdots \xi_n , n\geq 0\}$ 
 (where $T_0>0$ and $\{\xi_n\}$ are iid $\sim G$ and independent of $T_0$). In this
 case,  assume without loss of generality that $\bdist=\pp{0}$. 
However, if $\Gz>0$, 
any result not depending on $H$ must be restricted to fdds conditional
on the tail chain not having yet hit $\{0\}$. For example, consider
the trajectory of $(X_1,\dots,X_m)$, started from $X_0 = t$, through
the region $(t,\infty)^{m-2} \times [0,\delta] \times (t,\infty)$,
where $t$ is a high level. The tail chain would model this as a path
through $(0,\infty)^{m-2} \times \{0\} \times (0,\infty)$, which
requires specifying $H$ to control transitions away from $\{0\}$.

This raises the question of how to interpret the first hitting time of
$\{0\}$ for $\bT$ in terms of the original Markov chain $\bX$. Such
hitting times are important in the study of Markov chain point
process models of exceedance clusters based on the tail
chain. Intuitively, a transition to $\{0\}$ by $\bT$ represents a
transition from an extreme state to a non-extreme state by $\bX$. We
 make this notion precise in Section \ref{secextrbdry} by
viewing such transitions as downcrossings of a certain level we term
the ``extremal boundary.''

We assume  $\bX$ is a Markov chain on  $[0,\infty)$ with transition kernel $ K \in D(G)$,
%and
$K^*$ is  a tail kernel associated with $G$ with unspecified
boundary distribution $H$, and $\bT$ is a Markov chain on
$[0,\infty)$ with kernel  $K^*$.
The finite-dimensional distributions of $\bX$, conditional on $X_0=y$, are given by 
\[ \P_y \big[ (X_1,\, \dots,\, X_m) \in d\bx_m \big] = \mckern{y}{dx_1} \mckern{x}{dx_2} \dotsm \mckern{x_{m-1}}{dx_m}, \] and analogously for $\bT$.

%%%%%%%%%%%%%%%%%%%%%%%%%%%%%%%%%%%%%%%%%%%%%%%%%%%%%%%%%%%%

\subsection{FDDs Conditional on the Intial State}

\label{secfddvcsimple}

\newcommand{\fddt}[3]{\mckern[\pi^{(t)}_{#3}]{#1}{#2}}
\newcommand{\fdd}[3]{\mckern[\pi_{#3}]{#1}{#2}}
\newcommand{\fddvt}[3]{\mckern[\mu^{(t)}_{#3}]{#1}{#2}}
\newcommand{\fddv}[3]{\mckern[\mu_{#3}]{#1}{#2}}

Define the \sid{conditional distributions}
\begin{equation}\label{eqn:defPi}
\fddt{u}{\cdot}{m} = \P_{tu}\biggl[\bigg(\frac{X_1}{t},\dots,\frac{X_{m}}{t}\bigg) \in \cdot\, \biggr] \quad\text{and}\quad
\fdd{u}{\cdot}{m} = \P_u\big[( T_1,\dots, T_m ) \in \cdot\, \big],
\quad m\geq 1,
 \end{equation}
on $[0,\infty) \times \Bor[0,\infty]^m$. 
We \sid{consider when} $\pi^{(t)}_m \wc \pi_m$ on $[0,\infty]^m$
pointwise in $u$. \sid{If $\Gz=0$}, this is a direct consequence of
the domain of attraction condition \eqref{eqkerndoa}, but
if $\Gz>0$, more thought is required. 
We begin by
 restricting the convergence to
the smaller space $\E'_m := (0,\infty]^{m-1} \times
[0,\infty]$. Relatively compact sets in $\E'_m$ are contained in rectangles 
$[\vect{a},\binfty] \times [0,\infty]$, where $\vect{a} \in
(0,\infty)^{m-1}$.  

\begin{thm} \label{thmfddvcsimple}
Let $u_t = u(t)$ be a non-negative function such that $u_t \goesto u>0$ as $t\goesto\infty$. 
\begin{enumerate}
	\item[(a)] The restrictions to $\E'_m$, 
\begin{equation}\label{eqn:defMu}
 \fddvt{u}{\cdot\,}{m}: = \fddt{u}{\cdot\, \cap \E'_m}{m} \quad
 \text{and} \quad
\fddv{u}{\cdot\,}{m} := \fdd{u}{\cdot\, \cap \E'_m}{m},
\end{equation}
satisfy 
\begin{equation} \label{eqfddvcsimple}
\fddvt{u_t}{\cdot\,}{m} \vconv \fddv{u}{\cdot\,}{m} \mt{in} \mplus(\E'_m) \qquad(t\goesto\infty).
\end{equation}

	\item[(b)] 
	If $\Gz=0$, we have 
\begin{equation}\label{eqn:addone}
 \fddt{u_t}{\cdot\,}{m} \wc \fdd{u}{\cdot\,}{m}
 \mt{on} [0,\infty]^m \qquad(t\goesto\infty).
\end{equation}
\end{enumerate}
\end{thm}

\begin{pf}
The Markov structure suggests an induction argument facilitated by Lemma \ref{lemintconv} (p.~\pageref{lemintconv}). 
Consider (a) first. 
If $m=1$, then \eqref{eqfddvcsimple} above reduces to \eqref{eqkernunifconv}. 
Assume $m\geq 2$, and let $f\in \contfn (\E'_m)$. Writing $\E'_m = (0,\infty]\times \E'_{m-1}$, we can find $a>0$ and $B\in \cpt(\E'_{m-1})$ such that $f$ is supported on $[a,\infty]\times B$. 
Now, observe that
\ba 
\fddvt{u_t}{\cdot}{m}(f) &= \int_{(0,\infty]}  \mckern{tu_t}{tdx_1} \int_{\E'_{m-1}} \mckern{tx_1}{tdx_2} \dotsm \mckern{tx_{m-1}}{tdx_m}
 \; f(\bx_m) \\ 
 &= \int_{(0,\infty]}  \mckern{tu_t}{tdx_1} 
  \int_{\E'_{m-1}} \fddvt{x_1}{d(x_2,\ldots,x_m)}{m-1} \; f(\bx_m) . 
\ea 
Defining 
\ba
 h_t(v) = \int_{\E'_{m-1}} \fddvt{v}{d\bx_{m-1}}{m-1} \; f(v,\bx_{m-1}) \quad \text{and}\quad 
 h(v) = \int_{\E'_{m-1}} \fddv{v}{d\bx_{m-1}}{m-1} \; f(v,\bx_{m-1}) \,, 
\ea
 the previous expression becomes \[ \fddvt{u_t}{\cdot}{m}(f) 
 = \int_{\sid{(0,\infty]}}  \mckern{tu_t}{tdv} \; h_t(v).
 \]
Now, suppose $v_t \goesto v>0\,$: we verify 
\begin{equation} \label{eqhtgoestoh}
h_t(v_t) \conv h(v). 
\end{equation}
By continuity, we have $f(v_t,\bx_{m-1}^t) \goesto f(v,\bx_{m-1})$ whenever $\bx_{m-1}^t \goesto \bx_{m-1}$, and the induction hypothesis provides
 $\mckernl[\mu_{m-1}^{(t)}]{v_t}{\cdot} \vconv \mckernl[\mu_{m-1}]{v}{\cdot}$. 
Also, $f(x,\cdot)$ has compact support $B$ (without loss of generality, $\mckernl[\mu_{m-1}]{v}{\bdry B} = 0$).  
Combining these facts, \eqref{eqhtgoestoh} follows from Lemma \ref{lemintconv} (b). 
Next, since the $h_t$ and $h$ have common compact support $[a,\infty]$, and recalling from Propostion \ref{propkernconvunif} that 
$\mckernl{tu_t}{t\,\cdot} \wc \mckernl[K^*]{u}{\cdot}$,  
Lemma \ref{lemintconv} (a) yields
\[ \fddvt{u_t}{\cdot}{m}(f) \conv \int_{(0,\infty]}  \mckern[K^*]{u}{dv} \; h(v) = \fddv{u}{\cdot}{m}(f). \]

Implication (b) follows from essentially the same argument. For $m\geq 2$, suppose $f\in \cbfn[0,\infty]^m$. Replacing $\mu$ by $\pi$ and $\E'_{m-1}$ by $[0,\infty]^{m-1}$ in the definitions of $h_t$ and $h$, we have 
\[ \fddt{u_t}{\cdot}{m}(f) = 
\int_{[0,\infty]} \mckern{tu_t}{tdv} \; h_t(v). 
 \] This time 
Lemma \ref{lemintconv} (a) shows that $h_t(v_t) \goesto h(v)$ if
$v_t \goesto v>0$, and since $\mckernl[K^*]{u}{(0,\infty]}=1$, 
resorting to Lemma \ref{lemintconv} (a) once more yields 
\[ 
\fddt{u_t}{\cdot}{m}(f)
\conv \int_{[0,\infty]}  \mckern[K^*]{u}{dv} \; h(v) = \fdd{u}{\cdot}{m}(f). \qedhere \] 
\end{pf} 

If $\Gz>0$, then 
$\mckernl[K^*]{u}{(0,\infty]} = 1-\Gz < 1$,
\sid{and} for \eqref{eqn:addone} to hold 
would require knowing the behaviour of $h_t(v_t)$ when $v_t\goesto 0$
as well. 
Behaviour near zero is controlled by an asymptotic condition related to the boundary distribution $\bdist$. 
Previous work handled this using the regularity condition
discussed in Section \ref{secregcond}.

%%%%%%%%%%%%%%%%%%%%%%%%%%%%%%%%%%%%%%%%%%%%%%%%%%%%%%%%%%%%

\subsection{The Extremal Boundary}
\label{secextrbdry}

The normalization employed in the domain of attraction condition \eqref{eqkerndoa} suggests that, starting from a large state $t$, the extreme states are approximately scalar multiples of $t$. For example, we would consider a transition from $t$ into $(t/3,2t]$ to remain extreme. 
Thus, we think of states which can be made smaller than $t\delta$ for any $\delta$, if $t$ is large enough, as non-extreme. In this context, the set $[0,\sqrt{t}]$ would consist of non-extreme states. 

Under \eqref{eqkerndoa}, a tail chain path
through $(0,\infty)$ models the original chain $\bX$ travelling among
extreme states, and  
all of the non-extreme states are compacted into the state $\{0\}$ in the state space of $\bT$. Therefore, if $\bX$ is started from an extreme state, the portion of the tail chain depending solely on $G$ is informative up until the first time $\bX$ crosses down to a non-extreme state. If $\Gz=0$, such a transition would become more and more unlikely as the initial state increases in which case $G$ provides a complete description of the behaviour of $\bX$ in any finite number of steps following a visit to an extreme state (Theorem \ref{thmfddvcsimple} (b)). 

Drawing upon this interpretation, we develop a rigorous formulation of the distinction between extreme and non-extreme states, and we recast Theorem \ref{thmfddvcsimple} as convergence on the unrestricted space $[0,\infty]^m$ of the conditional fdds, given that $\bX$ has not yet reached a non-extreme state. 

\begin{defn}
Suppose 
$K\in D(G)$. An \emph{extremal boundary} for $K$ is a non-negative function $\ebt$ defined on $[0,\infty)$, satisfying $\lim{t\goesto\infty} \ebt = 0$ and 
\begin{equation}\label{eqextrbdryconv}
\mckern{t}{t\,[0,\ebt]} \conv \Gz \qquad\text{as }\ t\goesto\infty.
\end{equation}
\end{defn}
\noindent Such a function is guaranteed to exist by Lemma \ref{lemwcseq} (p.~\pageref{lemwcseq}). 

If $\Gz=0$, then $\ebt \equiv 0$ is a trivial choice. 
For any function $0\leq \ebt \to  0$, we have $\limsup_{t\to\infty}
\mckernl{t}{t\,[0,\ebt]} \leq \Gz$, so 
\eqref{eqextrbdryconv}  is equivalent to 
\begin{equation}\label{eqn:altbdryconv}
\liminf_{t\to\infty}\,\mckern{t}{t\,[0,\ebt]} \geq \Gz. \end{equation}
If $\ebt$ is an extremal boundary, it follows that any function $0\leq \tilde{\eb}(t) \to  0$ with $\tilde{\eb}(t) \geq \ebt$ for $t\geq t_0$ is also an extremal boundary for $K$.  Taking $\tilde{\eb}(t)=\vee_{s\geq t}\,\eb(s)$ shows that without loss of generality, we can assume $\ebt$ to be non-increasing. 

The extremal boundary has a natural formulation in terms of the update function. As in \eqref{equpdstdform}, let  $\psi(y, (Z,W)) = Zy + \phi(y,W) $ be an
update function in canonical form, where $y$ is extreme. If $Z>0$ then
the next state is approximately $Zy$, another extreme
state. Otherwise, if $Z=0$, the next state is $\phi(y,W)$, and
a transition from an extreme to a non-extreme state has taken place. This
suggests choosing  an extremal boundary whose order is between $t$ and
$\phi(t, \sid{w})$.

\begin{prop} \label{propupdextrbdry}
Suppose $\psi(y,(Z,W))$ is an update function in canonical form as in
\eqref{equpdstdform}.
If $\zeta(t)>0$ is a function on $[0,\infty)$ such that 
\begin{equation} \label{eqexbdupd}
{\phi(t,w)}/{\zeta(t)} \conv 0 
\end{equation}
 as $t\goesto\infty$ whenever $w\in B$ for which $\P[W\in B] = 1$,
 then 
$\liminf_{t\to \infty} \mckernl{t}{[0,\zeta(t)]} \geq \Gz.$ Provided
$\lim_{t\to\infty} \zeta(t)/t = 0$, an extremal boundary is given by $\ebt:=\zeta(t)/t$.
\end{prop}
Thus if $\phi(t,w)=o(\zeta (t))$ and $\zeta (t) =o(t)$ then $\zeta
(t)/t$ is an extremal boundary. \sid{For example, if $\psi(y,(Z,W)) =
  Zy + W$, so that $\phi(t,w)=w$, then choosing $\zeta (t)$ to be
  any function $\zeta(t)\to\infty$ such that $\zeta (t)=o(t)$    makes $\zeta (t)/t$ an extremal boundary. 
Choosing $\zeta(t) = \sqrt{t}$, we find that $\ebt =
  1/\sqrt{t}$ is an extremal boundary. } 

\begin{pf}
Since 
\ba
\P\big[ \psi(t) \leq \zeta(t) \,,\, Z=0 \big] &= \P\big[ \phi(t,W) \leq \zeta(t) \,,\, Z=0 \big] \geq \P\big[ \abs{\phi(t,W)} \leq \zeta(t) \,,\, Z=0 \big] \\ 
 &\geq \P\big[ Z=0 \big] - \P\bigg[ \frac{\abs{\phi(t,W)}}{\zeta(t)} > 1 \bigg] \conv \P\big[Z=0\big],  
\ea 
we have 
\[  \liminf_{t\goesto\infty} \,\mckern{t}{[0,\zeta(t)]} = \liminf_{t\goesto\infty} \,\P\big[ \psi(t) \leq \zeta(t) \big] \geq \P\big[Z=0 \big]. \qedhere 
\]
\end{pf}

We will need an extremal boundary for which
\eqref{eqextrbdryconv} still holds upon replacing the initial state
$t$ with $tu_t$, where $u_t \to u>0$.
Compare the following extension with Proposition \ref{propkernconvunif}. 

\begin{prop}
If $K\in D(G)$, then there exists an extremal boundary $\eb^*(t)$ such that 
\begin{equation} \label{eqextrbdryconvunif}
\mckern{tu_t}{t\,[0,\eb^*(t)]} \conv \Gz \qquad\text{as }\ t\goesto\infty
\end{equation}
for any non-negative function $u_t = u(t) \goesto u>0$. 
\end{prop}

\noindent We will refer to $\eb^*$ as a \emph{uniform extremal boundary}. 

\begin{pf}
Let $\ebt$ be an extremal boundary for $K$. 
As a first step, fix $u_0>1$, and suppose $\inv{u_0} < u <
u_0$. Define $\tilde{\eb}(t) = u_0\,\eb(t\inv{u_0})$. Now, if $u_t\goesto
u$, then $\eb_{\{u\}}(t) := u_t\,\eb(tu_t)$ satisfies
\eqref{eqextrbdryconvunif}, since \[
\mckern{tu_t}{t\,[0,\eb_{\{u\}}(t)]} = \mckern{tu_t}{tu_t\,[0,\eb(tu_t)]}
\conv \Gz. \]  
Here $\eb_{\{u\}}$ depends on the choice of function $u_t$. However,
since we eventually have $\inv{u_0}<u_t<u_0$ for $t$ large enough, it
follows that $\tilde{\eb}(t) > \eb_{\{u\}}(t)$ for such $t$. Hence,
$\tilde{\eb}(t)$ satisfies \eqref{eqextrbdryconvunif} for any
$u_t\goesto u$ with $\inv{u_0} < u < u_0$.  

Next, we remove the restriction in $u_0$ via a diagonalization
argument. For $k=2,3,\dots$, let $\eb_k(t)$ be extremal boundaries such
that $\mckernl{tu_t}{t\,[0,\eb_k(t)]} \goesto \Gz$ whenever $u_t\goesto
u$ for $u\in (\inv{k},k)$, and put $\eb_0 = \eb_1 = \eb$. Next, define the
sequence $\{(s_k,x_k) :\, k=0,1,\dots \}$ inductively as
follows. Setting $s_0 = 0$ and $x_0 = \eb_0(1)$, choose $s_k \geq
s_{k-1}+1$ such that $\eb_j(t) \leq \inv{k} \smin x_{k-1}$ for all
$j=0,\dots,k$ whenever $t\geq s_k$, and put $x_k = \max\{\eb_j(s_k) :\,
j=0,\ldots,k\}$. Note that $x_k \leq \inv{k} \smin x_{k-1}$, so $x_k
\downarrow 0$, and $s_k \uparrow \infty$. Finally, set \[ \eb^*(t) =
\sum_{k=0}^\infty x_k \, \indfn_{[s_k,\,s_{k+1})}(t). \] Observe that
$0 \leq \eb^*(t) \downarrow 0$, and suppose $u_t \goesto u>0$. Then
$u\in (\inv{k_0},k_0)$ for some $k_0$, so
$\mckernl{tu_t}{t\,[0,\eb_{k_0}(t)]} \goesto \Gz$, and for $k\geq k_0$,
our construction ensures that whenever $s_k \leq t < s_{k+1}$, we have
$\eb_{k_0}(t) \leq \eb_{k_0}(s_k) \leq x_k = \eb^*(t)$. Therefore, $\eb^*(t)
\geq \eb_{k_0}(t)$ for $t\geq s_{k_0}$, so $\eb^*$ satisfies
\eqref{eqextrbdryconvunif}.  
\end{pf}

Henceforth, we  assume  any $K\in D(G)$ is accompanied by a
uniform extremal boundary denoted by $\ebt$, and we consider 
extreme states on the order of $t$ to be $(t\ebt,\infty]$.  
If $\Gz=0$, then all positive states are extreme states.
We now use the extremal boundary  to reformulate the convergence of
\sid{Theorem} \ref{thmfddvcsimple} on the larger space $[0,\infty]^m$. Put $\E'_m(t) =
(\ebt,\infty]^{m-1}\times [0,\infty]$, so that $\E'_m(t) \uparrow
\E'_m = (0,\infty]^{m-1}\times [0,\infty]$.   Recall the notation $\mu_m^{(t)}$ and
$\mu_m^*$ from \eqref{eqn:defPi}, \eqref{eqn:defMu} in
\sid{Theorem} \ref{thmfddvcsimple} (p.~\pageref{thmfddvcsimple}).
 
\newcommand{\fddvyt}[3]{\mckern[\widetilde{\mu}^{(t)}_{#3}]{#1}{#2}}

\begin{thm} \label{thmfddvcyt}
Let $u_t = u(t)$ be a non-negative function such that $u_t \goesto
u>0$ as $t\goesto\infty$.  Taking   
\[ \fddvyt{u}{\cdot\,}{m} = \fddt{u}{\cdot\, \cap \E'_m(t)}{m}, \]
 we have 
\begin{equation*} %\label{eqfddvcsimple}
\fddvyt{u_t}{\cdot\,}{m} \vconv \fddv{u}{\cdot\,}{m} \mt{in} \mplus [0,\infty]^m \qquad (t\goesto\infty). 
\end{equation*}
%in as $t\goesto\infty$.
\end{thm}

\begin{pf}
Note that we can just as well write $\mckernl[\widetilde{\mu}_m^{(t)}]{u}{\cdot} = \mckernl[\mu_m^{(t)}]{u}{\cdot\, \cap \E'_m(t)}$. 
Suppose $m\geq 2$ and let $f\in \contfn[0,\infty]^m$. For $\delta>0$, define $A_\delta = (\delta,\infty]^{m-1}\times [0,\infty]$, and choose $\delta$ such that $\mckernl[\mu_m]{u}{\bdry A_\delta} = 0$. 
On the one hand, for large $t$ we have 
\ba
\fddvyt{u_t}{\cdot}{m}(f) &= 
\int_{[0,\infty]^m} f(\bx)\, \indfn_{\E'_m(t)}(\bx) \; \fddvt{u_t}{d\bx}{m}  
 \geq 
\int_{
\E'_m
} 
f(\bx)\, \indfn_{A_\delta}(\bx)\; \fddvt{u_t}{d\bx}{m} \\ 
 &\qquad\conv \int_{\E'_m} f(\bx)\, \indfn_{A_\delta}(\bx)\; \fddv{u}{d\bx}{m}
\ea 
as $t\goesto\infty$ by Lemma \ref{lemsubsetintconv} (p.~\pageref{lemsubsetintconv}). Letting $\delta \downarrow 0$ yields  
\begin{equation} \label{eqextrbdryliminf}
 \liminf_{t\goesto\infty}\; 
  \fddvyt{u_t}{\cdot}{m}(f)
 \,\geq \,
\fddv{u}{\cdot}{m}(f)
\end{equation}
 by monotone convergence. 
On the other hand, fixing $\delta$, we can decompose the space according to the first downcrossing of $\delta$: 
\begin{equation} \label{eqextrbdrydecomp}
 \fddvyt{u_t}{\cdot}{m}(f)
 = \int_{[0,\infty]^m} f(\bx)\, \indfn_{A_\delta}(\bx) \; \fddvyt{u_t}{d\bx}{m} + \sum_{k=1}^{m-1} \int_{[0,\infty]^m} f(\bx)\, \indfn_{A_\delta^k}(\bx)  \; \fddvyt{u_t}{d\bx}{m},  
\end{equation}
where $A_\delta^k = (\delta,\infty]^{k-1}\times [0,\delta]\times [0,\infty]^{m-k}$. On the subsets $A_\delta^k$ we appeal to the bound on $f$, say $M$, to obtain 
\[ \int_{[0,\infty]^m} f(\bx)\, \indfn_{A_\delta^k}(\bx)  \; \fddvyt{u_t}{d\bx}{m} \leq M\,\fddvyt{u_t}{A_\delta^k}{m}.  \]
Now, 
\begin{align} \label{eqextrbdryterm}
\fddvyt{u_t}{A_\delta^k}{m} &\leq \fddvt{u_t}{(\delta,\infty]^{k-1}\times(\ebt,\delta]}{k}  \\ \nonumber
 &\qquad\qquad = \fddvt{u_t}{(\delta,\infty]^{k-1}\times[0,\delta]}{k} - \fddvt{u_t}{(\delta,\infty]^{k-1}\times [0,\ebt]}{k}. 
\end{align}
Considering the second term, we have
\ba
&\fddvt{u_t}{(\delta,\infty]^{k-1}\times [0,\ebt]}{k} \\ 
 &= \int_{[0,\infty]} \mckern{tu_t}{tdx_1} \indfn_{(\delta,\infty]}(x_{1}) \cdots \int_{[0,\infty]} \mckern{tx_{k-2}}{tdx_{k-1}} \indfn_{(\delta,\infty]}(x_{k-1}) \; \mckern{tx_{k-1}}{t\,[0,\ebt]} \\
 &= \int_{\E'_{k-1}} \fddvt{u_t}{d\bx_{k-1}}{k-1} \; h_t(\bx_{k-1}), 
\ea 
where 
\[ h_t(\bx_{k-1}) = \mckern{tx_{k-1}}{t\,[0,\ebt]} \, \indfn_{(\delta,\infty]^{k-1}}(\bx_{k-1}). \]
Moreover, if $\bx_{k-1}^t \goesto \bx_{k-1} \in (\delta,\infty]^{k-1}$, then 
\[ h_t(\bx_{k-1}^t) = \mckern{tx_{k-1}^t}{t\,[0,\ebt]} \, \indfn_{(\delta,\infty]^{k-1}}(\bx_{k-1}^t) \conv \Gz\,\indfn_{(\delta,\infty]^{k-1}}(\bx_{k-1}), \] 
using the fact that $\ebt$ is a uniform extremal boundary. Since $\mckernl[\mu_{k-1}]{u}{\bdry (\delta,\infty]^{k-1}} = 0$ without loss of generality by choice of $\delta$, we conclude that 
\[ \fddvt{u_t}{(\delta,\infty]^{k-1}\times [0,\ebt]}{k} \conv \Gz \cdot \fddv{u}{(\delta,\infty]^{k-1}}{k-1} = \fddv{u}{(\delta,\infty]^{k-1}\times \{0\}}{k} \] as $t\goesto\infty$. 
Now, let us return to \eqref{eqextrbdryterm}. Given any $\epsilon>0$, by choosing $\delta$ small enough, we can make 
\ba 
\fddvt{u_t&}{(\delta,\infty]^{k-1}\times(\ebt,\delta]}{k} \conv \fddv{u}{(\delta,\infty]^{k-1}\times [0,\delta]}{k} - \fddv{u}{(\delta,\infty]^{k-1}\times\{0\}}{k} \\ 
 &\leq \fddv{u}{(0,\infty]^{k-1}\times [0,\delta]}{k}  - \fddv{u}{(\delta,\infty]^{k-1}\times\{0\}}{k} \\
  &< \fddv{u}{(0,\infty]^{k-1} \times \{0\}}{k} + \frac{\epsilon}{2}  - \bigg(\fddv{u}{(0,\infty]^{k-1} \times\{0\}}{k} - \frac{\epsilon}{2} \bigg) = \epsilon, 
\ea 
\ie
\begin{equation} \label{eqextrbdrytermeps}
 \limsup_{t\goesto\infty} \; \fddvyt{u_t}{A_\delta^k}{m} < \epsilon, 
\end{equation}
for $k=1,\dots,m-1$. 
Therefore, 
\eqref{eqextrbdrydecomp} implies that, given $\epsilon'>0$, 
\ba
\limsup_{t\goesto\infty} \;
 \fddvyt{u_t}{&\cdot}{m}(f)
 \leq \int_{[0,\infty]^m} f(\bx)\, \indfn_{A_\delta}(\bx) \; \fddv{u}{d\bx}{m} + 
M \sum_{k=1}^{m-1} \limsup_{t\goesto\infty} \ \fddvyt{u_t}{A_\delta^k}{m} 
\\ &
< 
  \fddv{u}{\cdot}{m}(f)
 + \epsilon'
\ea 
for small enough $\delta$. Combining this with \eqref{eqextrbdryliminf} yields the result. 
\end{pf}

%%%%%%%%%%%%%%%%%%%%%%%%%%%%%%%%%%%%%%%%%%%%%%%%%%%%%%%%

\subsection{The Extremal Component}

Having thus formalized  
the distinction between extreme and non-extreme states, we
return to the question of phrasing a fdd limit result
for $\bX$ when $\bdist$ is unspecified. The extremal boundary
allows us to interpret the first hitting time of $\{0\}$ by the tail
chain as 
approximating the time of the first transition  from extreme down
to non-extreme. In this terminology, Theorem \ref{thmfddvcyt} provides
a  result, given that such a transition has yet to occur.  

Define the first hitting time of a non-extreme state  
 \[ \tau(t) = \inf \big\{n \geq 0 \,:\, X_n \leq t\ebt \big\}\,. \] 
 For a Markov chain started from $tu_t$, where $u_t\goesto u>0$, we
 have $tu_t > \ebt$ for large $t$,  
so $\tau(t)$ is the first downcrossing of the extremal boundary. 

For the tail chain $\bT$, put $\tau^* = \inf\{n\geq 0 :\, T_n
= 0\}$. Given $T_0>0$, write $\tau^* = \inf\{n\geq 1 :\, \xi_n
= 0\}$, where $\{\xi_n\}\sim G$ are iid and independent of
$T_0$, 
\ie $\tau^*$ follows a Geometric distribution with
parameter $p=\Gz$. Thus, $\P[\tau^* = m] = p(1-p)^{m-1}$ for
$m\geq 1$  if $p>0$, and $\P[\tau^* = \infty] = 1$ if $p=0$.  
 Theorem \ref{thmfddvcyt} becomes 
\begin{equation} \label{eqfddconvtau}
\P_{tu_t} \bigl[t^{-1}{\bX_m} \in \cdot \,,\; \tau(t) \geq m \bigr]
\vconv  \P_u \big[ \bT_m \in \cdot \,,\; \tau^* \geq m \bigr],  
\end{equation}
implying that $\tau^*$ approximates $\tau(t)$:
\begin{equation}\label{corconvdctime}
 \P_{tu_t} \big[ \tau(t) \in \cdot\, \big] \wc \P \big[ \tau^* \in
 \cdot\, \big],\quad (t\goesto\infty,\ u_t \to u>0).
 \end{equation}

So if $\Gz>0$,  $\bX$  takes an average of approximately $\inv{\Gz}$
steps to return to a non-extreme state.  
but  if $\Gz=0$,  $\P_{tu_t}[\tau_1 \leq m] \to 0$ for any $m\geq 1$
so 
starting from a larger and larger initial state, it will take longer
and longer for $\bX$ to cross down to a non-extreme state.  

Let  $\bT^*$ be the tail chain associated with $(G,\pp{0})$. 
For $\{\xi_n\}\sim G$  iid and independent of $T^*_0$,
\begin{equation} \label{eqtailchainorig}
T_n^* = T_0^*\, \xi_1 \cdots \xi_n .
\end{equation}
We restate \eqref{eqfddconvtau} in terms of a process derived from
$\bX$, called the {\it extremal component\/} of $\bX$,
 whose fdds converge weakly to those of $\bT^*$. 
The extremal component is the part of $\bX$ whose asymptotic behavior
is controlled by $G$ alone.

\begin{defn}
The \emph{extremal component} of $\bX$ relative to $t$
is the process $\bX^{(t)}$ defined for $t>0$ as  
\[ X_n^{(t)} = X_n \cdot \ind{n < \tau(t)}\,, \qquad n=0,1,\dots  . \] 
\end{defn}

\noindent Observe that $\bX^{(t)}$ is a Markov chain on $[0,\infty)$ with transition kernel 
\[ \mckern[K^{(t)}]{x}{A} = \begin{cases} \mckern{x}{A\cap (t\ebt,\infty]} + \pp{0}(A)\cdot \mckern{x}{[0,t\ebt]} & x>t\ebt \\ \pp{0}(A) & x\leq t\ebt \end{cases}\,. \] It follows that $\mckernl[K^{(t)}]{t}{t\,\cdot} \wc G$ as $t\goesto\infty$, and additionally that $\mckernl[K^{(t)}]{t}{\{0\}} \goesto \Gz$.

The relation between the component processes $\bX^{(t)}$,  $\bT^*$
and the complete ones is
\[ \P_{tu_t} \big[ t^{-1}{\bX_m^{(t)}} \in \cdot \bgiven \tau(t) > m \big]
= 
\P_{tu_t} \big[ t^{-1}{\bX_m} \in \cdot \bgiven \tau(t) > m \big] \] 
%\qquad\text{and}\qquad  
and
\[ \P_u \big[\bT^*_m \in \cdot \bgiven \tau^*> m \big] = \P_u \big[\bT_m \in \cdot \bgiven \tau^*> m \big]. \]

\newcommand{\fddect}[3]{\mckern[\widetilde{\pi}^{(t)}_{#3}]{#1}{#2}}

\begin{thm} \label{thmfddextrcomp}
Let $u_t = u(t)\geq 0$ satisfy $u_t \goesto u>0$ as $t\goesto\infty$. Then on $[0,\infty]^m$, 
\[ \fddect{u_t}{\cdot\,}{m} := \P_{tu_t} \bigg[
 \bigg(\frac{X_1^{(t)}}{t},\ldots,\frac{X_m^{(t)}}{t}\bigg) \in \cdot\,
\bigg] \wc \P_u \big[ (T^*_1,\,\ldots,\, T^*_m ) \in \cdot\,
\big]
\qquad (t\goesto\infty). \]
\end{thm}

\begin{pf}
Suppose $m\geq 2$ and $f \in \cbfn [0,\infty]^m$, and assume first that $f\geq 0$. Then $f\in\contfn[0,\infty]^m$ as well, since the space is compact. 
Recall the notation of Theorem \ref{thmfddvcyt}.
Conditioning on $\tau(t)$, we can write  
\ba
\fddect{u_t}{&\cdot}{m}(f) = 
  \int_{(0,\infty]^{m}} f(\bx_m) \; \fddect{u_t}{d\bx_{m}}{m} + \sum_{k=1}^{m} \int_{(0,\infty]^{k-1} \times \{0\}^{m-k+1}} f(\bx_m) \; \fddect{u_t}{d\bx_m}{m} \\
 &= \int_{(0,\infty]^{m}} f(\bx_m) \; \fddect{u_t}{d\bx_{m}}{m} + \sum_{k=1}^{m} \int_{(0,\infty]^{k-1} \times \{0\}} f(\bx_{k},\,0,\,\dots,\,0) \; \fddect{u_t}{d\bx_k}{k}
\ea 
by the Markov property. 
Since 
\ba
 \fddect{u_t}{&\cdot\, \cap\, (0,\infty]^{m}}{m} 
= \P_{tu_t} \big[ t^{-1}{\bX^{(t)}_m} \in \cdot\;,\; \tau(t) > m \big] 
 = \P_{tu_t} \big[ t^{-1}{\bX_m} \in \cdot\, \cap (\ebt,\infty]^m \big] \\ 
 &= \fddvyt{u_t}{\cdot \times [0,\infty]}{m+1}\,,
\ea
the first term becomes 
\[ \fddvyt{u_t}{\cdot}{m+1}(f) \conv \fddv{u}{\cdot}{m+1}(f) = \int_{(0,\infty]^m} f(\bx_m) \; \fdd{u}{d\bx_m}{m} = \int_{(0,\infty]^m} f(\bx_m) \; \P_u\big[\bT^*_m\in d\bx_m \big] \]
as $t\goesto\infty$. 
Next, for any $A \cont [0,\infty]^k$ measurable, write $A_0 = \{\bx_{k-1} \,:\, (\bx_{k-1}\,,\,0) \in A\} \cont [0,\infty]^{k-1}$, and observe that 
\ba
\fddect{u_t}{ A \,\cap\, &(0,\infty]^{k-1}\times \{0\}}{k} = \P_{tu_t} \big[ t^{-1}{\bX_{k-1}^{(t)}} \in A_0 \cap (0,\infty]^{k-1}\,,\, X_{k}^{(t)} = 0 \big] \\ 
 &= \P_{tu_t} \big[ t^{-1}{\bX_{k-1}} \in A_0 \cap (\ebt,\infty]^{k-1}\;,\; {t}^{-1}{X_{k}} \leq \ebt \big] 
 \\ &
 = \fddvyt{u_t}{A_0 \times [0,\infty]}{k} - \fddvyt{u_t}{A_0\times [0,\infty]^2}{k+1}. 
\ea 
Applying this reasoning to the terms in the summation 
%and changing variables 
yields 
\ba
 \int_{[0,\infty]^{k}} & f(\bx_{k-1},\,0,\,\dots,\,0) \; \fddvyt{u_t}{d\bx_k}{k}  - \int_{[0,\infty]^{k+1}} f(\bx_{k-1},\,0,\,\dots,\,0) \; \fddvyt{u_t}{d\bx_{k+1}}{k+1} \\
 &\conv  \int_{[0,\infty]^{k}} f(\bx_{k-1},\,0,\,\dots,\,0) \; \fddv{u}{d\bx_k}{k}  - \int_{[0,\infty]^{k+1}} f(\bx_{k-1},\,0,\,\dots,\,0) \; \fddv{u}{d\bx_{k+1}}{k+1} \\
 &\quad= \int_{(0,\infty]^{k-1}\times \{0\}} f(\bx_k,\,0,\,\dots,\,0) \; \fdd{u}{d\bx_k}{k} =  \int_{(0,\infty]^{k-1} \times \{0\}^{m-k+1}} f(\bx_m) \; \P_u\big[ \bT^*_m \in d\bx_m \big]. 
\ea 
Combining these limits shows that 
$ \EP_{tu_t} f\big(   {t}^{-1}  {\bX^{(t)}_m} \big) \conv \EP_u f(\bT^*_m) $,
%\[ \fddect{u_t}{\cdot}{m}(f) \conv \fdd{u}{\cdot}{m}(f) \] 
as $t\goesto\infty$. 
Finally, if $f$ is not non-negative, then write $f=f_+ - f_-\,$. Since each of $f_+$ and $f_-$ is non-negative, bounded, and continuous, we can apply the above argument to each. 
\end{pf}

%%%%%%%%%%%%%%%%%%%%%%%%%%%%%%%%%%%%%%%%%%%%%%%%%%%%%%%

\section{The Regularity Condition}
\label{secregcond}

Previous work on the tail chain derives fdd convergence of $\bX$ to
$\bT^*$ under a single assumption analogous to our domain of
attraction condition \eqref{eqkerndoa}. As we observed in Section
\ref{secfddvcsimple},  when
$\Gz=0$, 
fdd convergence of $\{t^{-1}\bX\}$ 
follows directly, but when $\Gz>0$, it \sid{was} common to assume an additional
\sid{technical} condition which made
\eqref{eqkerndoa} imply fdd convergence to $\bT^*$  as
well.  This condition, which we refer to as the
``regularity condition,''  is an asymptotic convergence assumption prescribing the boundary distribution to be $\bdist=\pp{0}$. 
We consider equivalences between different forms 
appearing in the literature, in terms of both kernels and update functions, 
and show that, under the regularity condition,
 the extremal behaviour of $\bX$ is asymptotically the same as that of
 its extremal component $\bX^{(t)}$.  

In cases where $\Gz>0$, Perfekt \cite{perfekt1994extremal,perfekt1997extreme} requires 
that 
\begin{equation} \label{eqregcondperfekt}
\lim_{\delta\downarrow 0}\; \limsup_{t\goesto\infty} \sup_{u\in[0,\delta]}\mckern{tu}{(t,\infty]} = 0, 
\end{equation} 
while Segers \cite{segers2007multivariate} stipulates that the chosen update function corresponding to $K$ must be of at most linear order in the initial state: 
\begin{equation} \label{eqregcondsegers} 
\limsup_{t\goesto \infty} \, \sup_{0\leq y\leq t} \inv{t}\psi(y,v) <
\infty, \qquad (v\in B_0,\ \P[V\in B_0]=1).
\end{equation}
Smith \cite{smith1992extremal} used a variant of
\eqref{eqregcondperfekt}. 
We deem a formulation in terms of distributional convergence to be instructive in our context.  

\begin{defn}
A Markov transition kernel $K\in D(G)$ satisfies the \emph{regularity condition} if 
\begin{equation} \label{eqregcondkern}
\mckern{tu_t}{t\cdot\,} \wc \pp{0}(\cdot) 
\end{equation} 
on $[0,\infty]$ as $t\goesto\infty$ for any non-negative function $u_t = u(t) \goesto 0$. 
\end{defn}

Note that in \eqref{eqkernunifconv} (p.~\pageref{eqkernunifconv}), we had $u_t \to u>0$. 
We  interpret \eqref{eqregcondkern} as designating the boundary distribution $\bdist$ to be $\pp{0}$. 

We now consider the relationships between 
\eqref{eqregcondperfekt},
\eqref{eqregcondsegers}
and \eqref{eqregcondkern} ,
and propose an intuitive equivalent for update functions in canonical form. 

\begin{prop} \label{propregcond}
Suppose $K\in D(G)$, and let $\psi(\cdot,V)$ be an update function corresponding to $K$ such that 
\begin{equation} \label{equpdfnasreg}
\inv{t}\psi(t,v)\conv \xi(v) 
\end{equation}
whenever $v\in B$ for which $\P[V\in B]=1$, and $\xi\circ V\sim G$. 
 Then: 
\begin{enumerate}
	\item[(a)] Condition \eqref{eqregcondperfekt} is  necessary
          and sufficient 
for $K$ to satisfy the regularity condition \eqref{eqregcondkern}. 
	
	\item[(b)] Condition \eqref{eqregcondsegers} is sufficient
 for $K$ to satisfy the regularity condition \eqref{eqregcondkern}. 
	
	\item[(c)] If $\psi$ is in canonical form, \ie \[ \psi(y,(Z,W)) = Zy + \phi(y,W), \] then 
	$\psi$ satisfies \eqref{eqregcondsegers} if and only if 
$\phi(\cdot,w)$ is bounded on any neighbourhood of $0$ for each $w\in C$, a set for which $\P[W\in C] = 1$. 
\end{enumerate}
\end{prop}

\begin{pf}
(a)\ \ Assume \eqref{eqregcondperfekt}, and suppose $u_t\goesto 0$. We show $\mckernl{tu_t}{t(x,\infty]}\goesto 0$ for any $x>0$. Write \[ \omega(t,\delta) = \sup_{u\in[0,\delta]} \mckern{tu}{(t,\infty]}. \] 
Let $\epsilon>0$ be given, and
choose $\delta$ small enough that  $\limsup_{t\goesto\infty} \omega(t,\delta) < \epsilon/2$. 
Then for $t$ large enough that $u_t < \delta x$, we have 
\[ \mckern{tu_t}{t(x,\infty]} \leq \sup_{u\in[0,\delta x]} \mckern{tu}{t(x,\infty]} = \omega(tx, \delta) <  \limsup_{t\goesto\infty}\; \omega(t,\delta) + \epsilon/2 \]
for $t$ large enough. Our choice of $\delta$ implies that $\mckernl{tu_t}{t(x,\infty]} < \epsilon$. 

Conversely, assume that $K$ satisfies \eqref{eqregcondkern} but that \eqref{eqregcondperfekt} fails. Choose $\epsilon > 0$ and a sequence $\delta_n \downarrow 0$ such that $\limsup_{t\goesto\infty} 
\omega(t,\delta_n)
\geq \epsilon$ for 
$n=1,2,\dots$. 
Then for each $n$ we can find a sequence $t_k^n \goesto \infty$ as $k\goesto\infty$ such that $
\omega(t_k^n, \delta_n)
 \geq \epsilon$ for each $k$. 
 Diagonalize to find $k_1 < k_2 < \cdots$ such that $s_n = t_{k_n}^n \goesto\infty$ and  $
\omega(s_n,\delta_n)
  \geq \epsilon$ for all $n$. 
Finally, for $n=1,2,\dots$ choose $u_n \in [0,\delta_n]$ such that 
\[ \mckern{s_nu_n}{(s_n,\infty]} > \omega(s_n,\delta_n) - \epsilon/2, \]
and put $u(t) = \sum_n u_n \,\indfn_{[s_n,s_{n+1})}(t)$. Clearly $u(t) \goesto 0$, but $\mckernl{s_n u(s_n)}{(s_n,\infty]} \geq \epsilon/2$ for all $n$, contradicting \eqref{eqregcondkern}. 

(b)\ \ Write  $M(v) = \limsup_{t} \, \sup_{0\leq y\leq t}\, \inv{t}\psi(y,v)$. 
Since 
\[ \sup_{0\leq y\leq t} \inv{t}\psi(y,v) 
 = \sup_{0\leq y\leq \delta} \frac{\psi(t\inv{\delta}y,v)}{t\inv{\delta}}\,\inv{\delta} 
 \] for $\delta>0$, we have \[ \limsup_{t\goesto \infty} \sup_{0\leq y\leq \delta} \inv{t}\psi(ty,v)= \delta M(v). \] 
Now, 
 suppose $u_t \goesto 0$. 
 Given any $\delta>0$ we have 
 \[ \inv{t}\psi(tu_t,v) \leq \sup_{0\leq y\leq \delta} \inv{t}\psi(ty,v) \] provided $t$ is large enough, so $\limsup_t \inv{t}\psi(tu_t,v) \leq \delta M(v)$.
 Consequently, $\limsup_t \inv{t}\psi(tu_t,v) = 0$ for every $v$ such that $M(v)<\infty$. Under \eqref{eqregcondsegers}, this means that $\P\big[\inv{t}\psi(tu_t,V) \goesto 0 \big] = 1$, implying \eqref{eqregcondkern}.

(c)\ \ Suppose first that $\chi_w(a) = \sup_{0\leq y\leq a} \phi(y,w)<\infty$ for all $a>0$, whenever $w\in C$. Fixing $w\in C$ and $z\geq 0$, note that \[ \sup_{0\leq y\leq t} \inv{t}\psi(y,(z,w)) \leq z + \sup_{0\leq y\leq t} \inv{t}\phi(y,w), \] and observe for any $a>0$ that 
\[ \sup_{0\leq y \leq t} \inv{t} \phi(y,w) 
  \leq \Big(\sup_{0\leq y \leq a} \inv{t} \phi(y,w)  \Big) \smax \Big(\sup_{a\leq y \leq t}\inv{y} \phi(y,w)  \Big) 
   \leq \inv{t} \chi_w(a) \smax \Big( \sup_{a\leq y }\, \inv{y}\phi(y,w)  \Big). \] 
Choosing $a$ large enough that $\sup_{a\leq y }\, \inv{y}\phi(y,w) \leq 1$, say, it follows that 
\[ \limsup_{t\goesto\infty} \sup_{0\leq y\leq t} \inv{t}\psi(y,(z,w)) \leq z+1, \] so $v=(z,w) \in B_0$. Therefore $\P[(Z,W) \in B_0] \geq \P[Z\geq 0,\, W\in C] = 1$.

Conversely, suppose there is a set $D$ with $\P[W\in D]>0$ such that $w\in D$ implies $\chi_w(a)=\infty$ for some $0<a<\infty$. Since $\sup_{0\leq y\leq t} \inv{t}\psi(y,(z,w)) \geq \inv{t} \chi_w(t)$, we have $[0,\infty)\times D \cont \comp{B_0}$, contradicting \eqref{eqregcondsegers}. 
\end{pf}

The exclusion of necessity from part (b) 
results from the fact that a kernel $K$ does not uniquely specify an update function
$\psi$. Even when
 $K$ satisfies the regularity condition
\eqref{eqregcondkern}, it may be possible to choose
 a nasty update function $\psi$ which
satisfies \eqref{equpdfnasreg}, but not \eqref{eqregcondsegers}.
However, in such cases there may exist a
different update function $\psi'$ corresponding to $K$ which does satisfy \eqref{eqregcondsegers}.

Here is an example of such a situation. We exhibit an update function $\psi$ for which 
(i) \eqref{equpdfnasreg} holds;
(ii) \eqref{eqregcondsegers} fails because condition (c) in
Proposition \ref{propregcond} fails; but yet
(iii) the corresponding kernel satisfies the regularity condition 
\eqref{eqregcondkern}. Furthermore, we present a different choice of update function corresponding to the same kernel which satisfies \eqref{eqregcondsegers}.
Define $\psi(y,V=(Z,W)) = Zy + \phi(y,W)$, where \[ \phi(y,w) =
\sum_{k=1}^\infty k \cdot \ind{yw=1/k} \] and $W\sim U(0,1)$. (i) Since
$\phi(t,w) = 0$ for $t>1/w$, it is clear that $\psi$ satisfies
\eqref{equpdfnasreg} with $\xi = Z$.  (ii)
Observe that for any $w\in (0,1)$, $\phi(\cdot,w)$ is unbounded on the interval $[0,1]$. Therefore, by part (c) of Proposition \ref{propregcond}, \eqref{eqregcondsegers}  cannot hold for $\psi$. 
(iii) However, the corresponding kernel does satisfy the regularity
condition \eqref{eqregcondkern}. Suppose $u_t\goesto 0$ and 
$a>0$ is arbitrarily large.  Write  \[ \P\big[\inv{t} \psi(tu_t,(Z,W)) > x \big] =
\P\big[Z u_t + \inv{t} \phi(tu_t,W) > x \big] \leq \P\big[\inv{t}
\phi(tu_t,W) > x' \big] + \P[Z > a], 
\]
choosing $0<x' < x-au_t$.  
Since for any $t$, $ \{w\,:\,  \phi(tu_t,w) > tx'  \} \cont \{ \inv{(tu_tk)} \,:\, k=1,2,\dots \},
$ a set of measure 0 with respect to $\P[W\in \cdot\,]$,
\eqref{eqregcondkern} follows by letting $a\to\infty$.
On the other hand, the update function $\psi'(y,Z) = Zy$ does
 satisfy \eqref{eqregcondsegers}, and for any $y$, \[
 \P\big[\psi'(y,Z) \neq \psi(y,(Z,W)) \big] = \P\big[W \in
 \{\inv{(yk)} \,:\; k=1,2,\dots\} \big] = 0, \] so $\psi'$ does indeed
 correspond to $K$.  

\sid{The regularity condition \eqref{eqregcondkern}} restricts
attention to Markov chains for which the probability of returning to
an extreme state in the next $m$ steps after falling below the
extremal boundary is asymptotically negligible. 
For such chains, as well as those for which $\ebt\equiv 0$ is an extremal boundary for $K$, 
$\bX$ has the same asymptotic behaviour as its extremal component, as
described next.  

\begin{thm} \label{thmfddregcond} 
Suppose $\bX \sim K$ with $K\in D(G)$, and let $\rho$ be a metric on
$\R^m$. If $\ebt\equiv 0$ is an extremal boundary for $K$, or if $K$ satisfies the regularity condition
\eqref{eqregcondkern}, then for any $\epsilon>0$ we have  
\begin{equation} \label{eqextrcompslutsky}
 \P_{tu_t} \bigg[ \rho\Big(\frac{\bX_m^{(t)}}{t}\,,\, \frac{\bX_m}{t} \Big) >
 \epsilon \bigg] \conv 0 \qquad (t\to\infty,\ u_t \to u>0). 
\end{equation}
Consequently, 
\begin{equation} \label{eqfddconv}
\P_{tu_t} \bigg[ \bigg( \frac{X_1}{t},\dots,\frac{X_m}{t}\bigg) \in \cdot\, \bigg] \wc \P_u \big[ (T^*_1,\,\dots,\, T^*_m ) \in \cdot\, \big] \qquad (t\to\infty,\ u_t \to u>0). 
\end{equation}
\end{thm}

\noindent First let us extend the regularity condition to higher-order
transition kernels.

\begin{lem} \label{lemregcondhd}
If $K$ satisfies \eqref{eqregcondkern}, then so do the $m$-step transition kernels $K^m$. 
\end{lem}

\begin{pf}
This is established by induction.
Let $u_t\goesto 0$ and $f\in \cbfn[0,\infty]$. For $m\geq 2$, we have 
\[ \mckern[K^m]{tu_t}{\cdot}(f) = \int_{[0,\infty]} \mckern[K^{m-1}]{tu_t}{tdv} \int_{[0,\infty]}  \mckern{tv}{tdx} \; f(x) . \] 
Assume that $\mckernl[K^{m-1}]{tu_t}{t\,\cdot} \wc \pp{0}$; \eqref{eqregcondkern} implies that $\int \mckernl{tv_t}{tdx}\; f(x) \goesto f(0)$ whenever $v_t\goesto 0$. Therefore, by Lemma \ref{lemintconv} (a) (p.~\pageref{lemintconv}), we conclude that \[ \mckern[K^m]{tu_t}{\cdot}(f) \conv f(0) = \pp{0}(f). \qedhere
 \]
\end{pf}

\begin{pf}[Proof of Theorem \ref{thmfddregcond}]
Suppose $\epsilon>0$ and $u_t\to u>0$. Write 
\[ \P_{tu_t} \big[ \rho(t^{-1}{\bX_m^{(t)}}\,,\, t^{-1}{\bX_m}) > \epsilon \big] = \sum_{k=1}^m \, \P_{tu_t} \big[
\rho ( t^{-1}{\bX_m^{(t)}}\,,\, \inv{t}{\bX_m} ) > \epsilon
\;,\; \tau(t) = k \big]. \] 
Since $X_j = X_j^{(t)}$ while $j<\tau(t)$, for the $k$-th summand to converge to 0, it is sufficient that  
\[ \P_{tu_t} \big[ \lvert {X_j^{(t)}}/{t} - {X_j}/{t}\rvert >
\delta \;,\; \tau(t) = k \big] 
= \P_{tu_t}\big[ {X_j}/{t} > \delta \;,\; \tau(t) = k \big] \conv 0 \] 
for $j=k,\dots,m$ and any $\delta>0$. 
If $j=k$, we have 
\[ \P_{tu_t}\big[ {X_j}/{t} > \delta \;,\; \tau(t) = k \big] \leq
\P_{tu_t}\big[ 
{X_k}/{t} > \delta \,,\, {X_k}/{t} \leq \ebt \big] = 0 \] for large $t$. 
For $j>k$, recalling the notation of Theorem \ref{thmfddvcyt}, 
\ba
 \P_{tu_t}\big[ {X_j}/{t} &> \delta \;,\; \tau(t) = k \big] 
 = \int_{\E'_k(t)}  \indfn_{[0,\ebt]}(x_k)\,\P_{tu_t}\big[ {X_j}/{t} > \delta  \big| {\bX_k}/{t} = \bx_k \big] \; 
  \P_{tu_t}\big[{\bX_k}/{t}\in d\bx_k \big] \\
 &= \int_{[0,\infty]^k} \P_{tx_k}\big[ X_{j-k} > t\delta \big] \, \indfn_{[0,\ebt]}(x_k) \; \fddvyt{u_t}{d\bx_k}{k} 
\ea 
using the Markov property. We claim that this intergral $\goesto 0$ as $t\goesto\infty$. If $\ebt\equiv 0$, this follows directly. Otherwise, 
recall that $\mckernl[\widetilde{\mu}_k^{(t)}]{u_t}{\cdot} \vgoesto \mckernl[\mu_k]{u}{\cdot}$, 
and consider $h_t(\bx_k) = \P_{tx_k}[ X_{j-k} > t\delta ] \, \indfn_{[0,\ebt]}(x_k)$. 
Suppose $\bx^{(t)} \goesto \bx \in [0,\infty]^k$. If $x_k>0$, 
then $h_t(\bx^{(t)}) = 0$ for large $t$ because $\ebt\goesto 0$. 
Otherwise, if $x_k=0$, 
we have $h_t(\bx^{(t)}) \goesto 0$ since Lemma \ref{lemregcondhd} implies that  
$\P_{tx_k^{(t)}}[ X_{j-k} > t\delta ] \goesto 0$
as $t\goesto \infty$. Lemma \ref{lemintconv} (b) 
establishes \eqref{eqextrcompslutsky}; \eqref{eqfddconv} follows by Slutsky's theorem. 
\end{pf}

Therefore, $\bX$ converges to $\bT^*$ in fdds under (a) $\Gz=0$, (b)
$\Gz>0$ combined with \eqref{eqregcondkern}, or (c) $\Gz>0$ combined
with the extremal boundary $\ebt \equiv 0$. 
In either case, we will be able to replace the extremal component $\bX^{(t)}$ with the complete chain $\bX$ in the results of Sections \ref{secextvconv} and \ref{secjmrv}. 
However, that $\ebt \equiv 0$ is an extremal boundary, and consequently that \eqref{eqfddconv} holds, does not imply
the regularity condition to hold, regardless of $\Gz$; in particular,
a kernel for which $\Gz=0$ need not satisfy
\eqref{eqregcondkern}. This is illustrated in Example
\ref{exunifcounter}.

%%%%%%%%%%%%%%%%%%%%%%%%%%%%%%%%%%%%%%%%%%%%%%%%%%%%%%%%%%

\section{Convergence of the Unconditional FDDs}

\subsection{Effect of a Regularly Varying Initial Distribution} \label{secextvconv}

So far our  convergence results required that 
the initial state become large, and the only distributional
assumption  was that the transition kernel $K$ determining
$\bX$ be attracted to some distribution $G$.  
To obtain a result for the unconditional distribution of
$(X_0,\dots,X_m)$, we require an additional assumption about 
how likely the initial observation $X_0$ is to be large.
Using Lemma \ref{lemextvconv}, 
the results of the previous sections extend  to
multivariate regular variation on the cone $\E_m = (0,\infty] \times
[0,\infty]^m$ when the distribution of $X_0$ has a regularly varying
tail. 
This cone is smaller than the cone  $[0,\infty]^{m+1} \backslash
\{\bz\}$ traditionally
employed in extreme value theory, because the kernel domain of attraction condition
\eqref{eqkerndoa} is uninformative when the initial state is not
extreme. This is analogous to the setting of the Conditional Extreme
Value Model considered in \cite{das2011conditioning,
  heffernan2007limit}.  

\begin{prop} \label{propextvconv}
Assume $\bX \sim K$ with $K\in D(G)$, and
$X_0 \sim H$, where $H$ is a distribution on $[0,\infty)$ with a regularly
varying tail. 
This means that as $t \to \infty$, for some scaling function $b(t)\to\infty$,
\[ tH\big(b(t)\cdot\big) \vconv
\nu_\alpha(\cdot) \quad\text{\ in\ \ } \mplus(0,\infty], \] 
 where
$\nu_\alpha(x,\infty]=x^{-\alpha}$ and $\alpha>0$. Define the measure
$\nu^*$ on  
$\E_m = (0,\infty]\times [0,\infty]^m$ 
%$\E_m$ 
by 
\begin{equation} \label{eqlimmeas} 
 \nu^* \big(dx_0, d\bx_m\big) = \nu_\alpha(dx_0)\,\P_{\sid{x_0}}\big[ (T^*_1,\,\dots,\, T^*_m) \in d\bx_m \big]. 
\end{equation}
Then, for $m=1,2,\ldots$, the following convergences take place as $t\goesto\infty$: 
\begin{enumerate}
\item[(a)] In $\mplus ((0,\infty]^m \times [0,\infty])$, 
\[
 t\, \P \big[   b(t)^{-1}
 ({X_0}\,,\,{X_1}\,,\,\ldots\,,\,{X_{m}}) \in \cdot\, \cap
 %(0,\infty] \times \E'_m  
 (0,\infty]^m \times [0,\infty]
 \big]  \vconv
\nu^*\big(\cdot\,\cap\, 
(0,\infty]^m \times [0,\infty]
\big).
\]

\item[(b)] In $\mplus(\E_m)$,
\[
 t\, \P
 \big[{b(t)}^{-1} ({X_0^{(b(t))}}\,,\, {X_1^{(b(t))}}\,,\,\ldots\,,\, {X_{m}^{(b(t))}})
 \in 
\cdot\, \big]  \vconv 
\nu^*(\cdot).
\]

\item[(c)] 
If either $\Gz=0$, $\ebt\equiv 0$ is an extremal boundary, or $K$ satisfies the regularity condition
\eqref{eqregcondkern}, then in $\mplus(\E_m)$,
\[
 t\, \P \big[ {b(t)}^{-1} ({X_0}\,,\,{X_1}\,,\, \ldots\,,\,{X_{m}}) \in \cdot\,
 \big]  \vconv 
\nu^*(\cdot).
\]

\item[(d)] In $\mplus(0,\infty]$,
\[
t\, \P \big[{X_0}/{b(t)} \in dx_0 \;,\; \tau(b(t)) \geq m \big] \vconv \big(1-\Gz\big)^{m-1} \cdot \nu_\alpha(dx_0) 
.\]
\end{enumerate}
\end{prop}

\begin{rmk}
These convergence statements may be reformulated equivalently as, say, 
\[
\P \big[{b(t)}^{-1}( {X_0},{X_1},\ldots,{X_{m}}) \in \cdot \bgiven X_0 > b(t) \big]  \wc 
\P\big[ (T^*_0,\, T^*_1,\, \ldots,\, T^*_m ) \in \cdot\, \big],
\]
where $T^*_0\sim$ Pareto$(\alpha)$. This is the form considered by Segers \cite{segers2007multivariate}. 
\end{rmk}

\begin{pf}
Apply Lemma \ref{lemextvconv} (p.~\pageref{lemextvconv}) to the results of Theorems \ref{thmfddvcsimple},  \ref{thmfddextrcomp} and \ref{thmfddregcond}, and \eqref{corconvdctime}. 
\end{pf}

In the case $m=1$, $\E_1$ is a rotated version of $\E_\sqcap$ used in
the conditional extreme value model in \cite{das2011conditioning,das2011detecting} and 
the limit can be expressed as \[ \nu^*\big((x_0,\infty]\times [0,x_1] \big)  = 
\sid{\int_{x_0}^\infty \nu_\alpha (du)P[\xi\leq x_1/u]} = 
x_0^{-\alpha}\, \P\big[\xi \leq x_1/x_0 \big] - x_1^{-\alpha}\, \EP
\xi^{\alpha} \ind{\xi \leq x_1/x_0}  \] for $(x_0,x_1) \in (0,\infty]
\times [0,\infty]$, where $\xi\sim G$ (with $\EP
\xi^{\alpha} \leq \infty$). Since
\[
\nu^*\big((x_0,\infty]\times\{0\}\big) = x_0^{-\alpha}\, \P[\xi=0]
\qquad\text{and}\qquad \nu^*\big((0,\infty]\times (x_1,\infty]\big) =
x_1^{-\alpha}\,\EP \xi^\alpha,
\] 
 sets on the
$x_0$-axis incur mass proportional to $\Gz$, and sets bounded away
from this axis are weighted accordng to $\EP\xi^{\alpha}$.
A consequence of the second observation is that  
\[ \liminf_{t\to\infty} \, t\,\P\big[{X_1}/{b(t)} > x \big] \geq \EP \xi^\alpha \cdot x^{-\alpha}. \] 
Thus, knowledge concerning the tail behaviour of
$X_1$  imposes a restriction on the  distributions $G$ to which $K$
can be attracted via the $\alpha$-th moment. For example, if
$t\,\P[X_1/b(t) \in \cdot] \vgoesto \nu_\alpha$, then we must have
$\EP\xi^\alpha \leq 1$; this property will be examined further in the
next section 
and appears in various forms in Segers \cite{segers2007multivariate} and Basrak and Segers \cite{basrak2009regularly}, in the stationary setting.

%%%%%%%%%%%%%%%%%%%%%%%%%%%%%%%%%%%%%%%%%%%%%%%%%%%%%%%

\subsection{Joint Tail Convergence}
\label{secjmrv}
What additional assumptions
are necessary for convergences (b) and (c) of the previous result
to take place on the larger cone $\E^*_m = [0,\infty]^{m+1} \backslash
\{\bz\}$?  
This was considered by Segers \cite{basrak2009regularly,
  segers2007multivariate}
for stationary Markov chains. 
In (b), the
dependence on the extremal threshold and hence on $t$ means we are in the context of a triangular array and not, strictly
speaking, in the setting of joint regular variation.
However, the result is still  useful, for example, to derive a point process convergence via the Poisson transform \cite[p.~183]{resnick2007heavy}.

As a first step, we characterize convergence on the larger cone by
decomposing it into smaller, more familiar cones. This is similar to Theorem
6.1 in \cite{segers2007multivariate} and one of the implications of
Theorem 2.1 in \cite{basrak2009regularly}.  \sid{As a convention in what follows, set
$[0,\infty]^0\times A=A$}. Also, recall the notation $\E_m = (0,\infty]\times [0,\infty]^m$. 

\begin{prop} \label{propcharmrv}
Suppose $\rvect{Y}_t = (Y_{t,0}\,,\, Y_{t,1}\,,\,\ldots,\, Y_{t,m})$ is a random vector on $[0,\infty]^{m+1}$ for each $t>0$. Then there exists a non-null Radon measure $\mu^*$  on $\E^*_m = [0,\infty]^{m+1} \backslash \{\bz\}$ such that 
\begin{equation} \label{eqpropmrv}
 t\,\P\big[ ( Y_{t,0}\,,\, Y_{t,1}\,,\, \dots,\, Y_{t,m} ) \in
 \cdot\, \big] \vconv \mu^*(\cdot) \mt{in} \mplus(\E^*_m)
\qquad (t\to \infty)
\end{equation}
	if and only if
	for $j=0,\ldots,m$ there exist Radon measures $\mu_j$ on $\E_j = (0,\infty] \times [0,\infty]^{j}$, not all null, such that 
\begin{equation} \label{eqpropsubmrv}
 t\,\P\big[ (Y_{t,j}\,,\, \dots,\, Y_{t,m} ) \in \cdot\, \big] \vconv \mu_{m-j}(\cdot) \mt{in} \mplus (\E_{m-j}). 
\end{equation}
The relation between the limit measures is the following: 
\[ \mu_{m-j}(\cdot) = \mu^*\big([0,\infty]^{j} \times \cdot\,\big) \mtt{on} \E_{m-j} \] for $j=0,\ldots,m$, 
and 
\[ \mu^*\big(\comp{[\bz,\bx]}\big) = \sum_{j=0}^m \mu_{m-j}\big( (x_j,\infty]\times [0,x_{j+1}] \times \cdots \times [0,x_m] \big) \mtt{for} \bx\in \E^*_m. \] 
Furthermore, given $j \in \{0,\dots,m-1\}$, if $A\cont [0,\infty]^{m-j} \backslash \{0\}^{m-j}$ is relatively compact, then $\mu_{m-j}((0,\infty]\times A) < \infty$. 
\end{prop}

\begin{pf} 
Assume first that \eqref{eqpropmrv} holds. Fixing $j \in \{0,\dots,m\}$, define $\mu_{m-j}(\cdot) 
:= \mu^*([0,\infty]^{j} \times \cdot\, )$ (\ie $\mu_m = \mu^*$). Let $A\cont \E_{m-j}$ be relatively compact with $\mu_{m-j}(\bdry A) = 0$. Then $A^* = [0,\infty]^j \times A$ is relatively compact in $\E^*_m$, and $\bdry_{\E^*_m} A^* = [0,\infty]^j \times \bdry_{\E_{m-j}} A$, so $\mu^*(\bdry_{\E^*_m} A^*) = \mu_{m-j}(\bdry A) = 0$. Therefore, 
\[ t\,\P\big[ (Y_{t,j}\,,\, \dots,\, Y_{t,m} ) \in A \big] = t\,\P\big[ ( Y_{t,0}\,,\, \dots,\, Y_{t,m} ) \in A^* \big] \conv \mu^*(A^*) = \mu_{m-j}(A), \] establishing \eqref{eqpropsubmrv}. 

Conversely, suppose we have \eqref{eqpropsubmrv} for $j=0,\dots,m$. For $\bx \in (0,\infty]^{m+1}$, define 
\[ h(\bx) = \sum_{j=0}^m \mu_{m-j}\big( (x_j,\infty]\times [0,x_{j+1}] \times \cdots \times [0,x_m] \big). \]
Decompose $\comp{[\bz,\bx]}$ as a disjoint union
\begin{equation} \label{eqpfdecompx}
 \comp{[\bz,\bx]} = \bigcup_{j=0}^m \ [0,\infty]^{j} \times (x_j,\infty] \times [0,x_{j+1}] \times \cdots \times [0,x_m]\,, 
 \end{equation}
and  observe that at points of continuity of the limit,
\begin{equation} \label{eqpfdecomp}
 t\,\P\big[ \rvect{Y}_t \in \comp{[\bz,\bx]}\, \big] 
 = \sum_{j=0}^m \, t\,\P\big[ ( Y_{t,j}\,,\, \dots,\, Y_{t,m}
 ) \in (x_j,\infty] \times [0,x_{j+1}] \times \cdots \times
 [0,x_m] \, \big] \conv h(\bx) .
 \end{equation}
Hence, \eqref{eqpropmrv} holds with the limit measure $\mu^*$ defined
by $ \mu^*\big( \comp{[\bz,\bx]} \big) = h(\bx). $
Indeed, given $f\in \contfn(\E^*_m)$ we can find $\delta>0$ such that $\bx_\delta = (\delta,\dots,\delta)$ is a continuity point of $h$ and $f$ is supported on $\comp{[\bz,\bx_\delta]}$. 
Therefore, 
\[ t\,\EP f(\rvect{Y}_t) \leq \sup_{\bx\in\E^*_m} f(\bx) \,\cdot\, \sup_{t>0}\, t\,\P\big[\rvect{Y}_t \in \comp{[\bz,\bx_\delta]} \big] < \infty,  \] 
 implying that the set $\big\{ t\,\P[\rvect{Y}_t \in \cdot\, ] \,;\; t>0 \big\}$ is relatively compact in $\mplus(\E^*_m)$. Furthermore, if $t_k\,\P[\rvect{Y}_{t_k} \in \cdot\, ] \goesto \mu$ and $s_k\,\P[\rvect{Y}_{s_k} \in \cdot\, ] \goesto \mu'$ as $k\goesto \infty$, then $\mu=\mu'=\mu^*$ on sets $\comp{[\bz,\bx]}$  which are continuity sets of $\mu^*$ by \eqref{eqpfdecomp}. This extends to measurable rectangles in $\E^*_m$ bounded away from $\bz$ whose vertices are continuity points of $h$, leading us to the conclusion that $\mu=\mu'=\mu^*$ on $\E^*_m$. 
 
Moreover, since we can decompose $\comp{[\bz,\bx]}$ for any $\bx\in\E^*_m$ as in \eqref{eqpfdecompx}, it is clear that $\mu^*$ is non-null iff not all of the $\mu_j$ are null. 

Finally, for $1\leq j\leq m-1$, if $A\cont [0,\infty]^{m-j} \backslash \{0\}^{m-j}$ is relatively compact, then it is contained in $\comp{[(0,\dots,0),(x_{j+1},\dots,x_m)]}$ for some $(x_{j+1},\dots,x_m) \in (0,\infty]^{m-j}$. Applying \eqref{eqpfdecompx} once again, we find that 
\ba
 \mu_{m-j}\big((0,\infty] &\times A\big) = \mu^*\big([0,\infty]^j\times (0,\infty] \times A \big) \\ 
 &\leq \sum_{k=j+1}^m \mu^*\big([0,\infty]^{j+1} \times [0,\infty]^{k-j-1} \times (x_k,\infty] \times [0,x_{k+1}] \times \cdots \times [0,x_m] \big) \\ 
  &= \sum_{k=j+1}^m \mu_{m-k}\big((x_k,\infty] \times [0,x_{k+1}] \times \cdots \times [0,x_m] \big) < \infty.  \qedhere 
\ea  
\end{pf}

Consequently, the extension of the convergences in Proposition \ref{propextvconv} to the larger cone $\E_m^*$ follows from regular variation of the marginal tails. 

\begin{thm} \label{thmjmrv}
Suppose $\bX \sim K\in D(G)$, and let $b(t)\goesto\infty$ be a scaling function and $\alpha>0$. Then 
\begin{equation} \label{eqpropmrvx}
t\, \P \big[
{b(t)}^{-1}({X_0^{(b(t))}},\,{X_1^{(b(t))}},\,\ldots,\,{X_{m}^{(b(t))}})
\in \cdot \, \big]  \vconv \mu^*(\cdot)  \mt{in} \mplus(\E^*_m) \qquad (t\goesto\infty), 
\end{equation}
where 
\[ \mu^*\big\rvert_{\E_m}(dx_0,d\bx) = \nu_\alpha(dx_0)\,\P_{x_0}\big[
(T^*_1,\,\dots,\, T^*_m) \in d\bx_m \big] = \nu^*(dx_0, d\bx), \] 
 if and only if 
\begin{equation} \label{eqpropmargrv}
 t\,\P\big[{X_j^{(b(t))}}/{b(t)} \in \cdot\, \big] \vconv c_j\, \nu_\alpha(\cdot) 
\end{equation}
 in $\mplus(0,\infty]$, with $c_0=1$ and $(\EP\xi^{\,\alpha})^j \leq c_j < \infty$ for $j=1,\dots,m$. \end{thm}

\begin{pf}
Assume first that \eqref{eqpropmrvx} holds.
It follows that 
\[ t\,\P\big[ X_0 > b(t)x] \conv \nu^*((x,\infty]\times [0,\infty]^m)
= x^{-\alpha} \] for $x>0$. Hence, $b(t) \in \RV_{1/\alpha}$, \sid{so by \eqref{eqpropmrvx} again, we have for}
$j\geq 1$ 
\[ 
 t\,\P\big[{X_j^{(b(t))}}> {b(t)} x \big] \conv \mu^*([0,\infty]^j\times (x,\infty] \times [0,\infty]^{m-j}) = c_j \, x^{-\alpha}, 
\]
and 
\ba
 c_j \geq \mu^*((0,\infty]&\times [0,\infty]^{j-1} \times (1,\infty] \times [0,\infty]^{m-j}) = \int_{(0,\infty]} \nu_\alpha(du) \; \P\big[\xi_1 \cdots \xi_j > \inv{u}] \\ 
 &= \EP (\xi_1 \cdots \xi_j)^{\alpha} = (\EP\xi^{\alpha})^j.  
 \ea 
Conversely, suppose that \eqref{eqpropmargrv} holds for
$j=0,\dots,m$. Lemma \ref{lemextvconv} implies that   in $\mplus(\E_{m-j})$,
\ba
 t\, \P \big[{b(t)}^{-1}({X_j^{(b(t))}} ,\ldots,{X_{m}^{(b(t))}}) \in (dx_0,d\bx) \, \big] 
 \vconv c_j\, &\nu_\alpha(dx_0)\, \P_{x_0}\big[ (T^*_1,\dots,T^*_{m-j}) \in d\bx \big] \\
&=: \mu_{m-j}\bigl((dx_0,d\bx)\bigr) 
\ea 
by the Markov property, and  Proposition \ref{propcharmrv} yields \eqref{eqpropmrvx}, with $\mu^* \rvert_{\E_m}(\cdot) = \mu_m(\cdot) = \nu^*(\cdot)$. 
\end{pf}

At the end of Section
\ref{secregcond}, cases were outlined in which we could replace $X_j^{(b(t))}$
by $X_j$.  Theorem \ref{thmjmrv} is most striking for these since it shows that
 for a Markov chain whose kernel is in a domain of
attraction, to obtain joint regular variation of the fdds it is enough
to know that the marginal tails are regularly varying.  
In particular, if $\bX$ has a regularly varying stationary
distribution then the fdds are jointly regularly varying. This result
was presented by Segers \cite{segers2007multivariate}, and Basrak and
Segers \cite{basrak2009regularly} showed that for a general stationary
process, joint regular variation of fdds is equivalent to the
existence of a ``tail process'' which reduces to the tail chain in the
case of Markov chains. However, what Proposition \ref{thmjmrv} emphasizes
is that it is the marginal tail behaviour alone, rather than
stationarity, which provides the link with joint regular variation.  

Theorem \ref{thmjmrv} also  extends the
observation made in Section \ref{secextvconv} 
that knowledge of the marginal tail behaviour for a Markov chain whose
kernel is in a domain of attraction constrains the class of possible
limit distributions $G$ via its moments. If a
particular choice of regularly varying initial distribution leads to  
$t\,\P[X_j > b(t)\,\cdot\,] \vgoesto a_j \nu_\alpha(\cdot)$, then we
have $\EP \xi^\alpha \leq a_j^{1/j}$. In particular, if $\bX$
admits a stationary distribution whose tail is $\RV_{-\alpha}$, then
$\EP \xi^\alpha \leq 1$.

%%%%%%%%%%%%%%%%%%%%%%%%%%%%%%%%%%%%%%%%%%%%%%%%%%%%%%%%%%

\section{Examples}
Our first example  illustrates the main results. 
\begin{exnum}
Let $V=(Z,\,W)$ be any random vector on $[0,\infty) \times \R$. 
Consider the update function $ \psi(y,V) = (Zy + W)_+ $ and its
canonical form 
$$\psi(y,V)=
Zy+\phi(y,W)= Zy + \big(W\,\ind{W > -Zy} - Zy\,\ind{W \leq -Zy}\big). 
$$
For $y>0$, the transition kernel has the form 
$\mckernl{y}{(x,\infty)} = \P\left[Zy+W > x\right]$. 
Since $\inv{t} \psi(t,V) = (Z + \inv{t}W)_+ \goesto Z$ a.s., we have $K \in D(G)$ with $G=\P[Z \in \cdot\, ]$. 
Furthermore, using Proposition \ref{propupdextrbdry}, the function $\gamma(t) \equiv \sqrt{t}$ is of larger order than $\phi(t,w)$, so $\ebt = 1/\sqrt{t}$ is an extremal boundary. 
Since $\phi(\cdot, w)$ is bounded on neighbourhoods of 0,  Proposition
\ref{propregcond} (c) implies
$K$  satisfies the regularity condition \eqref{eqregcondkern}.
Consequently, from Theorem \ref{thmfddregcond},
we obtain fdd convergence of $\sid{t^{-1}}\bX$ to $\bT^*$ as in
\eqref{eqfddconv}.
\end{exnum}

If $K$ does not satisfy the regularity condition
\eqref{eqregcondkern},
Theorem \ref{thmfddregcond} may fail to hold and 
starting from $tu$,  $\sid{t^{-1}}\bX$
may  fail to converge to $\bT^*$ started from $u$.

\begin{exnum} \label{exnonunif}
Let $V = (Z,\,W,\,W')$ be any non-degenerate random vector on $[0,\infty)^3$, and consider the Markov chain determined by the update function 
\[ \psi(y,V) = Zy + W\,\inv{y}\,\ind{y>0}+ W'\,\ind{y=0}. \] 
For $y>0$, the transition kernel is
$\mckernl{y}{(x,\infty)} = \P[Zy + W\inv{y} > x] $
and since $\inv{t} \psi(t,V) = Z + W t^{-2} \goesto Z$ a.s., we have $K\in D(G)$ with $G=\P[Z\in \cdot\,]$. 
Furthermore, using Proposition \ref{propupdextrbdry}, the function $\gamma(t) \equiv 1$ is of larger order than $\phi(t,w)$, so $\ebt = 1/t$ is an extremal boundary.

However, note that $\phi(y,(W,W')) = W \inv{y}\ind{y>0} + W'\ind{y=0}$ is unbounded near 0, implying that Segers' boundedness condition \eqref{eqregcondsegers} does not hold. In fact, our form of the regularity condition \eqref{eqregcondkern} fails for $K$.
Indeed, 
\[
 \mckern{tu_t}{t(x,\infty)} =\P[Ztu_t
+W/(tu_t)>tx]=\P[Zu_t+W/(t^2u_t) >x].
\]
Choosing $u_t = t^{-2}$ yields 
$ \mckernl{tu_t}{t(x,\infty)} \to \P[W>x]$. For appropriate $x$, this shows
\eqref{eqregcondkern} fails.

Not only does \eqref{eqregcondkern} fail but so does Theorem
\ref{thmfddregcond}, 
since the
asymptotic behaviour of $\bX$ is not the same as that of
$\bX^{(t)}$.
We show directly  that the conditional fdds of $\sid{t^{-1}}\bX$ fail
to converge to \sid{those of}
$\bT^*$. The idea is that if $X_k < \ebt = \inv{t}$, there is a
positive probability that $X_{k+1} > t$.  
We illustrate this for $m=2$. Take $f\in \cbfn [0,\infty]^2$ and
$u>0$. Observe if $X_0=tu>0$, from the definition of $\psi$,
$X_1=Z_1tu +W_1/(tu)$ and $X_2=Z_2X_1+(W_2/X_1)\ind{X_1>0}
+W'\ind{X_1=0}.$
Furthermore, on $\{Z_1>0\}$, we have $X_1>0$ and $X_2=Z_2X_1+W_2/X_1.$
On $\{Z_1=0,W_1>0\}$, $X_1>0$ and $X_2=Z_2X_1 +W_2/X_1$. 
On $\{Z_1=0,W_1=0\}$, we have $X_1=0$ and $X_2=W'$. Therefore
\ba
\EP_{tu} f( {X_1/t},\, {X_2/t} ) 
 &= \EP_{tu} f( {X_1/t},\, {X_2/t} )  \,\ind{ Z_1>0 }
+\EP_{tu} f( {X_1/t},\, {X_2/t} )  \,\ind{Z_1=0, W_1>0 } \\ 
 &\qquad\qquad + \EP_{tu} f( {X_1/t},\, {X_2/t} ) \, \ind{ Z_1=0, W_1=0} =A+B+C.
\ea 
For $A$, as $t\to\infty$, we have 
\ba
A &= \EP f\big(Z_1u+W_1/(t^2u),\, Z_2[Z_1u+W_1/(t^2u)] +W_2/[Z_1t^2u+W_1u^{-1}] \big) \,\ind{Z_1>0} \\ 
 &\qquad\conv \EP f(Z_1u,\,Z_1Z_2u) \,\ind{Z_1>0} ,
\ea
while for $B$ we obtain for $t\to\infty$,
\[ B = \EP f(W_1/t^2u,\, Z_2W_1/(t^2u) +W_2u/W_1) \, 
\ind{ Z_1=0, W_1>0 }  
\conv 
\EP  f(0,\, uW_2/W_1)\, 
\ind{ Z_1=0, W_1>0 } . \]
Finally for $C$,
\[ C = \EP f(0,\, W'_2/t) \, 
\ind{ Z_1=0, W_1=0 }
=\P[ Z_1=0, W_1=0 ] \,\EP f(0,\,W'_2/t) 
\conv \P[ Z_1=0, W_1=0 ] \, f(0,\,0). \]
Observe that $\lim_{t\to\infty}[A+B+C]\neq \EP_u f(T_1^*,T_2^*)=\EP f(uZ_1,uZ_1Z_2).$
\end{exnum}

In the final example, the \sid{conditional distributions of
  $t^{-1}\bX$} converge to those of the 
tail chain $\bT^*$, even though the regularity condition does not
hold. This includes cases for which $\Gz=0$ and $\Gz > 0$ with
extremal boundary $\ebt \equiv 0$.  

\begin{exnum}\label{exunifcounter} 
Let $\{(\xi_j,\eta_j),\, j\geq 1\}$ be iid copies of the non-degenerate random vector $(\xi,\eta)$ on $[0,\infty)^2$. Taking $V = (\xi,\,\eta)$, consider a Markov chain which transitions according to the update
function 
\[
 \psi(y,V) = \xi (y + {y}^{-1})\,\ind{y>0} + \eta\,\ind{y=0} =
 \xi\, y + \big(\xi\, \inv{y}\,\ind{y>0} + \eta\,\ind{y=0}\big) ,\] 
 where the last expression is the canonical form.
For $y>0$, the transition kernel is 
 \[ \mckern{y}{[0,x]} = \P\left[\xi (y + \inv{y}) \leq x\right] =
 \P\left[ \xi \leq {x}/(y+\inv{y}) \right].\] 
For $t>0$,  $\inv{t} \psi(t,V) = \xi(1+t^{-2}) \goesto \xi$ a.s., so $K\in D(G)$ with $G=\P[\xi\in \cdot\,]$. 
Note that $\phi(y,V) = \xi  \inv{y}\,\ind{y>0} + \eta\,\ind{y=0}$ is
unbounded near 0, implying that Segers' boundedness condition
\eqref{eqregcondsegers} does not hold. Also, our 
regularity condition \eqref{eqregcondkern} fails for $K$.  To see
this, write  
\[ 
\mckern{tu_t}{t(x,\infty)} = \P\left[\xi > {x}/{(u_t + \inv{(t^2u_t)})}\right]. \] 
Fix $x$ so that $\P[\xi>x]>0$ and choose $u_t = t^{-2}$. This yields $u_t + \inv{(t^2u_t)} = 1 + t^{-2}$,
implying that 
 \[
 \mckern{tu_t}{t(x,\infty)} = \P\left[\xi > {x}/{(1 +
     t^{-2})}\right] \geq \P[\xi > x]>0, \] 
so \eqref{eqregcondkern} fails for $K$.
However, since $\mckernl{t}{\{0\}} = \P[\xi=0] = \Gz$, the choice $\ebt \equiv
0$ satisfies the definition of an extremal boundary
\eqref{eqextrbdryconv}, even if $\Gz>0$. 
\sid{This leads to} fdd convergence of $\P_{tu}[t^{-1}\bX \in
\cdot\,]$ to $\P_u[\bT^* \in \cdot \,]$, and thus we learn that the
conclusion \eqref{eqfddconv}   of Theorem \ref{thmfddregcond}
may hold without \eqref{eqregcondkern} being true.

 We \sid{prove the fdd convergence} for $m=2$. 
For $u>0$, and $X_0=tu$, we have $X_1=\xi_1(tu+(tu)^{-1})$ and
$X_2=\xi_2(X_1+X_1^{-1})\,\ind{X_1>0} +\eta_2\,\ind{X_1=0}$. On $\{X_1>0\} = \{\xi_1>0\}$ we have $X_2=\xi_2(X_1+X_1^{-1})$. 
On $\{X_1=0\} = \{\xi_1=0\}$, we have $X_2=\eta_2$. Thus, as $t\to\infty$,
\ba
\EP_{tu} f(X_1/t,\,  X_2/t) \,\ind{X_1>0} &= \EP_{tu}  f\big(\xi_1[u+(t^2u)^{-1}],\,
   \xi_2 \xi_1 [u +(t^2u)^{-2}] 
+ \xi_2/ (\xi_1[t^2u +1/u]) \big) \, \ind{X_1>0} \\
 &\quad\conv \EP f(\xi_1u,\, \xi_1 \xi_2 u)\, \ind{\xi_1>0} ,
\ea 
while
\[ \EP_{tu} f(X_1/t,\,X_2/t) \,\ind{X_1=0} = \EP 
f(0,\,\eta_2/t)\, \ind{\xi_1=0} \conv \P[\xi_1=0]\, f(0,0). \]
We conclude that 
\[ \EP_{tu}  f(X_1/t,\, X_2/t ) \conv 
\EP f(u\xi_1,\, u\xi_1 \xi_2 ) =\EP_u f(T^*_1,\,T^*_2). \qedhere \]
\end{exnum}

%%%%%%%%%%%%%%%%%%%%%%%%%%%%%%%%%%%%%%%%%%

\section{Concluding Remarks}

We have thus placed the traditional tail chain model for the extremes of a Markov chain in a more general context through the introduction of the boundary distribution $\bdist$ as well as the extremal boundary. A common application of the tail chain model is in deriving the weak limits of exceedance point processes for $\bX$ \cite{basrak2009regularly, perfekt1994extremal, rootzén1988maxima}. 
%This would be a natural continuation of our development.
We will shortly use our results to develop a detailed description of the clustering properties of extremes of Markov chains by means of such point processes. 
Furthermore, as we have not employed stationarity in our finite-dimensional results, we propose to substitute the inherent regenerative structure of a Harris recurrent Markov chain for the traditional assumption of stationarity. 
%when considering point processes. 
Also, it would be interesting to explore the implications of choices of $\bdist$ other than $\pp{0}$.

%%%%%%%%%%%%%%%%%%%%%%%%%%%%%%%%%%%%%%%%%%%%%%%%%%%%%

\section{Appendix: \sid{Technical} Lemmas}
\label{seclem}

This section collects lemmas needed to prove
convergence of
integrals of the form $\int f_n \; d\mu_n$, assuming that $f_n \goesto
f$ and $\mu_n \goesto \mu$ in their respective spaces. An example is
the {\it second continuous mapping theorem\/} 
\cite[Theorem 5.5, p.~34]{billingsley1968convergence}.

\begin{lem} \label{lemunifcontmapping}
Assume $\E$ and $\E'$ are complete separable (cs) metric spaces,  and for
$n\geq 0$, $h_n :\E\goesto \E'$ are measurable.
 Put $ A = \{x \in \E : h_n(x_n) \goesto h_0(x) \text{\ whenever\ } x_n \goesto x \}. $
If $P_n$, $n \geq 0$ are probability measures on $\E$ with $P_n \wc P_0$, 
and  $h_n\to h_0$ almost uniformly in the sense that $P(A) = 1$, then  
$ P_n\circ\inv{h_n} \wc P_0\circ\inv{h_0}$   {in} $\E'.$
\end{lem}

\noindent The result provides a way to handle the convergence of a family of integrals.

\begin{lem} \label{lemintconv}
In addition to the assumptions of Lemma \ref{lemunifcontmapping},
require  $\E' = \R$ and $\{h_n,\, n\geq 0 \}$ is uniformly bounded, so that
$\sup_{n\geq 0}\sup_{x\in\E}|h_n(x)| < \infty$. 
\begin{enumerate}
	\item[\emph{(a)}] We have \[ \int_\E h_n \; dP_n \conv \int_\E h_0 \; dP_0. \]
	
	\item[\emph{(b)}] Suppose additionally that $\E$ is locally compact
          with a countable base (lccb), and $\mu_n \vgoesto \mu_0$ in
          $\mplus(\E)$ with $\mu_0 (\comp{A}) = 0$. If there exists a
          compact set $B\in\cpt(\E)$ with $\mu_0(\bdry B)=0$ such that
          $h_n(x)=0$, $n\geq 0$ whenever $x\not\in B$
          (\ie $B$ is a common compact support of each $h_n$),
          then \[ \int_{\E} h_n\;d\mu_n \conv \int_{\E} h_0\;d\mu_0. \]  
\end{enumerate}
\end{lem}

\begin{pf}
(a)\ \ 
 If $X_n\sim P_n$ for $n\geq 0$, then $h_n(X_n)\wc h_0(X_0)$. The uniform boundedness of the $h_n$ guarantees that $\EP h_n(X_n)\goesto \EP h_0(X_0)$. 

(b)\ \ View $B$ as a compact subspace of $\E$ inheriting the relative
topology.  Then, assuming $\mu(B) > 0$ to rule out a trivial case,
define probabilities on $B$ by
$ P_n(\cdot) = {\mu_n(\,\cdot\cap B)}/{\mu_n(B)},\,n\geq 0$.
Since $\mu_n(\,\cdot \cap B) \vconv   \mu_0(\,\cdot \cap B)$  by
Proposition 3.3 in \cite{feigin1996parameter}, and $B$
is compact, we get $P_n \wc P_0$.  
Denote by $h'_n$, $n\geq 0$, the  restriction of  $h_n$ to $B$. Observe that for any $x\in A\cap B$, we have 
$h'_n(x_n)\goesto h'(x)$ whenever $x_n \goesto x$ in $B$, and $P(\comp{A} \cap B) \leq \mu(\comp{A})/\mu(B) = 0$. 
Therefore, apply part (a) to obtain 
\[ \int_{\E} h_n\;d\mu_n = \int_\E h_n\, \indfn_B \;d\mu_n = \mu_n(B)\int_B h'_n\;dP_n \conv \mu_0(B)\int_{B} h'_0 \; dP_0 = \int_{\E} h_0 \; d\mu_0. \qedhere 
\] 
\end{pf} 

\noindent A convenient specialization of Lemma \ref{lemintconv} (b) is the following. 

\begin{lem} \label{lemsubsetintconv}
Suppose $\E$ is lccb and $\mu_n \vgoesto \mu$ in $\mplus(\E)$. If $f: \E \goesto \R$ is continuous and bounded, and $B \in \E$ is relatively compact with $\mu(\bdry B) = 0$, then \[ \int_B f\; d\mu_n \conv \int_B f\; d\mu. \]
\end{lem}

Take $h_n = f\indfn_B$ for  $n\geq 0$. Since $f\indfn_B$ is continuous except possibly on $\bdry B$, we have $\mu(\comp{A}) \leq \mu(\bdry B) = 0$.

The next result is used to extend convergence of substochastic
transition functions to multivariate regular variation on a larger
space.  

\begin{lem} \label{lemextvconv}
Let $\E\cont [0,\infty]^m$ and $\E' \cont [0,\infty]^{m'}$ be two nice (lccb) spaces. 
Suppose for $t\geq 0$ that
$\{\mckernl[p^{(t)}]{\,\cdot}{\cdot\,}\}_{t\geq 0}$, are substochastic
transition functions on $\E\times \Bor(\E')$. This means 
$\mckernl[p^{(t)}]{\cdot}{B}$ is a measurable function for any fixed $B \in
\Bor(\E')$, $\mckernl[p^{(t)}]{x}{\cdot}$ is a 
%substochastic probability 
measure for any $x\in\E$, 
and $\sup_{t\geq 0}\,\sup_{u\in \E}\, \mckernl[p^{(t)}]{u}{\E'} \leq 1. $
Assume there is a set $A\cont \E$ such that
\[ \mckern[p^{(t)}]{u_t}{\cdot\,} \vconv
\mckern[p^{(0)}]{u}{\cdot\,} \mt{in} \mplus(\E') \qquad (t\to\infty) \] 
 whenever $u_t\goesto u$ in $\E$ and $u \in A$.
Suppose also that $\{\nu^{(t)}\}_{t\geq 0}$ are measures on $\E$ such that $\nu^{(0)}(\comp{A}) = 0$,
and $ \nu^{(t)} \stackrel{v}{\to}\nu^{(0)} $ in $\mplus(\E)$.
Then, defining measures  $\mu^{(t)}$ for $t\geq 0$ on $\E\times \E'$ as 
 \[ \mu^{(t)}\big( du,dx \big) = \nu^{(t)}(du)\mckern[p^{(t)}]{u}{dx}\,, \] we have \[  \mu^{(t)} \vconv \mu^{(0)} \mt{in} \mplus (\E\times \E') \qquad (t\goesto\infty).  \]
\end{lem}

\begin{pf}
Let $f\in\contfn (\E\times \E')$; without loss of generality assume $f$ is supported on $K\times K'$, where $K\in \cpt(\E)$ and $K'\in \cpt(\E')$.
We have 
\[ \int_{\E\times\E'} \mu^{(t)}(du,dx) \; f(u,x) = \int_{\E} \nu^{(t)}(du)  \int_{\E'}  \mckern[p^{(t)}]{u}{dx} \; f(u,x). \]
For $t\geq 0$, write 
\[ \varphi_t(u) = \int_{\E'}  \mckern[p^{(t)}]{u}{dx} \; f(u,x)\]
and suppose $u_t \goesto u_0$ with $u_0\in A$; we verify that $\varphi_t(u_t)\goesto \varphi_0(u_0)$. 
Writing $g_t(x) = f(u_t,x)$, $t \geq 0,$  we have $g_t(x_t)\goesto
g_0(x_0)$ whenever $x_t\goesto x_0 \in \E'$ by the continuity of
$f$. Also, the $g_t$ are uniformly bounded by the bound on $f$ and
$g_t(x) = 0$ for all $t$ whenever $x \notin K'$. Furthermore,
without loss of generality we can assume that $\mckernl[p^{(0)}]{u}{\bdry K'}
= 0$. Now apply Lemma
\ref{lemintconv} (b) to obtain \[ \varphi_t(u_t) = \int_{\E'}
\mckern[p^{(t)}]{u_t}{dx} \; g_t(x) \conv \int_{\E'} \mckern[p^{(0)}]{u}{dx}
\; g_0(x) = \varphi_0(u). \] 
Since the $p^{(t)}$ are substochastic, and   $\varphi_t(u) = 0$ for all $t$ whenever $u\notin K$,
 the $\varphi_t$ are uniformly bounded by the bound on $f$ . Assume similarly that $\nu(\bdry K)=0$, and recall that $\nu(\comp{A})=0$. 
Apply Lemma \ref{lemintconv} (b) once more to conclude as $t \to \infty$ that
\[ \int_{\E\times\E'} \mu^{(t)}(du,dx) \; f(u,x) = \int_{\E}
\nu^{(t)}(du) \; \varphi_t(u)  
\conv  \int_{\E} \nu^{(0)}(du) \; \varphi_0(u) =  \int_{\E\times\E'} \mu^{(0)}(du,x) \; f(u,x) . \qedhere\]
\end{pf}

\noindent We conclude this section with a result used to verify the  existence of the extremal boundary.

\begin{lem} \label{lemwcseq}
Suppose $P_t$, $t\geq 0$ are probability measures on a cs metric space $\E$
such that $P_t \wc P_0$, and let $A\cont \E$ be measurable. Then there
exists a sequence of sets $A_t \downarrow \cl{A}$ such that $P_t(A_t)
\goesto P_0(\cl{A})$.  
\end{lem}

\begin{rmk}
Note that if $P(\bdry A) = 0$ then we can take $A_t = \cl{A}$. In the case of distribution functions $F_t \wc F$ on $\R^m$, taking $A = (-\infty,\bx]$ and metric $\rho=\rho_\infty$ shows that for any $\bx\in \R^m$ there exists $\bx_t \downarrow \bx$ such that $F_t(\bx_t) \goesto F(\bx)$. 
\end{rmk}

\begin{pf}
Let $\rho$ be a metric on $\E$, and consider sets $A_\delta = \{x \,:\, \rho(x,A) \leq \delta\}$. Recall that  $P_0(\bdry A_\delta) = 0$ for all but a countable number of choices of $\delta$, since $F(\delta) = P_0(A_\delta) - P_0(\cl{A})$ is a distribution function. First choose $\{\delta_k \,:\, k=1,2,\dots\}$ such that $0 < \delta_{k+1} \leq \delta_k \smin 1/(k+1)$ and $P_0(\bdry A_{\delta_k}) = 0$ for all $k$. Next, let $s_0 = 0$ and take $s_{k} \geq s_{k-1} + 1$, $k=1,2,\dots$ such that $P_t(A_{\delta_k}) > P_0(\cl{A}) - 1/k$ whenever $t\geq s_k$; this is possible since $P_t(A_{\delta_k}) \goesto P_0(A_{\delta_k}) \geq P_0(\cl{A})$ for all $k$. 
Finally, for $t>0$ set \[ A(t) = A_{\delta_1} \,\indfn_{(0,s_{1})}(t) + \sum_{k=1}^\infty A_{\delta_k} \,\indfn_{[s_k,s_{k+1})}(t). \] 
We claim that $A(t) \downarrow \cl{A}$ and that $P_t(A(t)) \goesto P_0(\cl{A})$ as $t\goesto\infty$. It is clear that $A(t) \supset A(t')$ for $t \leq t'$, and $\cap_t\, A(t) = \cap_k\, A_{\delta_k} = \cl{A}$. On the one hand, for large $t$ we have $A(t) \cont A_{\delta_k}$ for any $k$, so  
\[ \limsup_{t\goesto\infty}\, P_t(A(t)) \leq \limsup_{t\goesto\infty}\, P_t(A_{\delta_k}) \leq P_0(A_{\delta_k}). \] Letting $k\goesto\infty$ shows that $\limsup_t P_t(A(t)) \leq P_0(\cl{A})$. 
On the other hand, if $k(t)$ denotes the value of $k$ for which $s_k \leq t < s_{k+1}$, then 
\[ P_t(A(t)) = P_t(A_{\delta_{k(t)}}) > P_0(\cl{A}) - 1/k(t), \] so $\liminf_t P_t(A(t)) \geq P_0(\cl{A})$. 
Combining these two inequalities shows that $P_t(A(t)) \goesto P_0(\cl{A})$. 
\end{pf}

%%%%%%%%%%%%%%%%%%%%%%%%%%%%%%%%%%%%%%%%%%%%%%%%%%%

%% Bibliography

\nocite{bortot2003extremes}
\nocite{bortot1998models}
\nocite{davis1995point}
\nocite{de2006extreme}
\nocite{hsing1989extreme}
\nocite{kallenberg1997prob}
\nocite{resnick2007extreme}
\nocite{smith1997markov}
\nocite{yun1998extremal}
\nocite{yun2000distributions}


\begin{thebibliography}{10}

\bibitem{basrak2009regularly}
B.~Basrak and J.~Segers.
\newblock {Regularly varying multivariate time series}.
\newblock {\em Stochastic Processes and their Applications}, 119(4):1055--1080,
  2009.

\bibitem{billingsley1968convergence}
P.~Billingsley.
\newblock {\em Convergence of probability measures}.
\newblock Wiley, 1968.

\bibitem{billingsley:1971}
P.~Billingsley.
\newblock {\em Weak Convergence of Measures: {A}pplications in Probability}.
\newblock Society for Industrial and Applied Mathematics, Philadelphia, Pa.,
  1971.
\newblock Conference Board of the Mathematical Sciences Regional Conference
  Series in Applied Mathematics, No. 5.

\bibitem{billingsley1999convergence}
P.~Billingsley.
\newblock {\em Convergence of probability measures}.
\newblock Wiley-Interscience, 2nd edition, 1999.

\bibitem{bortot2000sufficiency}
P.~Bortot and S.~Coles.
\newblock {A sufficiency property arising from the characterization of extremes
  of Markov chains}.
\newblock {\em Bernoulli}, 6(1):183--190, 2000.

\bibitem{bortot2003extremes}
P.~Bortot and S.~Coles.
\newblock {Extremes of Markov chains with tail switching potential}.
\newblock {\em Journal of the Royal Statistical Society. Series B (Statistical
  Methodology)}, 65(4):851--867, 2003.

\bibitem{bortot1998models}
P.~Bortot and J.A. Tawn.
\newblock {Models for the extremes of Markov chains}.
\newblock {\em Biometrika}, 85(4):851, 1998.

\bibitem{das2011conditioning}
B.~Das and S.I. Resnick.
\newblock Conditioning on an extreme component: Model consistency with regular
  variation on cones.
\newblock {\em Bernoulli}, 17(1):226--252, 2011.

\bibitem{das2011detecting}
B.~Das and S.I. Resnick.
\newblock Detecting a conditional extreme value model.
\newblock {\em Extremes}, 14:29--61, 2011.

\bibitem{davis1995point}
R.A. Davis and T.~Hsing.
\newblock {Point process and partial sum convergence for weakly dependent
  random variables with infinite variance}.
\newblock {\em The Annals of Probability}, 23(2):879--917, 1995.

\bibitem{de2006extreme}
L.~de~Haan and A.~Ferreira.
\newblock {\em {Extreme value theory: an introduction}}.
\newblock Springer Verlag, 2006.

\bibitem{feigin1996parameter}
P.D. Feigin, M.F. Kratz, and S.I. Resnick.
\newblock Parameter estimation for moving averages with positive innovations.
\newblock {\em The Annals of Applied Probability}, pages 1157--1190, 1996.

\bibitem{heffernan2007limit}
J.E. Heffernan and S.I. Resnick.
\newblock {Limit laws for random vectors with an extreme component}.
\newblock {\em Annals of Applied Probability}, 17(2):537--571, 2007.

\bibitem{heffernan2004conditional}
J.E. Heffernan and J.A. Tawn.
\newblock {A conditional approach for multivariate extreme values}.
\newblock {\em Journal of the Royal Statistical Society. Series B (Statistical
  Methodology)}, 66(3):497--546, 2004.

\bibitem{hsing1989extreme}
T.~Hsing.
\newblock {Extreme value theory for multivariate stationary sequences}.
\newblock {\em Journal of Multivariate Analysis}, 29(2):274--291, 1989.

\bibitem{kallenberg1997prob}
O.~Kallenberg.
\newblock {\em {Foundations of modern probability}}.
\newblock Springer Verlag, 1997.

\bibitem{leadbetter1983extremes}
M.R. Leadbetter, G.~Lindgren, and H.~Rootz{\'e}n.
\newblock {\em Extremes and related properties of random sequences and
  processes}.
\newblock Springer series in statistics. Springer-Verlag, 1983.

\bibitem{perfekt1994extremal}
R.~Perfekt.
\newblock {Extremal behaviour of stationary Markov chains with applications}.
\newblock {\em The Annals of Applied Probability}, 4(2):529--548, 1994.

\bibitem{perfekt1997extreme}
R.~Perfekt.
\newblock {Extreme value theory for a class of Markov chains with values in
  $\R^d$}.
\newblock {\em Advances in Applied Probability}, 29(1):138--164, 1997.

\bibitem{resnick2007extreme}
S.I. Resnick.
\newblock {\em {Extreme values, regular variation, and point processes}}.
\newblock Springer Verlag, 2007.

\bibitem{resnick2007heavy}
S.I. Resnick.
\newblock {\em {Heavy-tail phenomena: probabilistic and statistical modeling}}.
\newblock Springer Verlag, 2007.

\bibitem{rootzén1988maxima}
H.~Rootz{\'e}n.
\newblock {Maxima and exceedances of stationary Markov chains.}
\newblock {\em ADV. APPL. PROB.}, 20(2):371--390, 1988.

\bibitem{segers2007multivariate}
J.~Segers.
\newblock {Multivariate regular variation of heavy-tailed Markov chains}.
\newblock {\em Arxiv preprint math/0701411}, 2007.

\bibitem{smith1992extremal}
R.L. Smith.
\newblock {The extremal index for a Markov chain}.
\newblock {\em Journal of applied probability}, 29(1):37--45, 1992.

\bibitem{smith1997markov}
R.L. Smith, J.A. Tawn, and S.G. Coles.
\newblock {Markov chain models for threshold exceedances}.
\newblock {\em Biometrika}, 84(2):249, 1997.

\bibitem{yun1998extremal}
S.~Yun.
\newblock {The extremal index of a higher-order stationary Markov chain}.
\newblock {\em Annals of Applied Probability}, 8(2):408--437, 1998.

\bibitem{yun2000distributions}
S.~Yun.
\newblock {The distributions of cluster functionals of extreme events in a
  dth-order Markov chain}.
\newblock {\em Journal of Applied Probability}, 37(1):29--44, 2000.

\end{thebibliography}
\end{document}